\documentclass[12 pt]{amsart}

\newtheorem{theorem}{Theorem}[section]

\newtheorem{proposition}[theorem]{Proposition}
\newtheorem{corollary}[theorem]{Corollary}
\newtheorem{problem}[theorem]{Problem}
\newtheorem{conjecture}[theorem]{Conjecture}

\theoremstyle{definition}

\newtheorem{example}[theorem]{Example}

\newtheorem{algorithm}[theorem]{Algorithm}

\theoremstyle{remark}

\usepackage{amscd,amssymb}

\begin{document}
\baselineskip 15 pt

\title[Coordinates and Automorphisms]
{Coordinates and Automorphisms of Polynomial
and Free Associative Algebras of Rank Three}

\author[Vesselin Drensky and Jie-Tai Yu]
{Vesselin Drensky and Jie-Tai Yu}
\address{Institute of Mathematics and Informatics,
Bulgarian Academy of Sciences,
1113 Sofia, Bulgaria}
\email{drensky@math.bas.bg}
\address{Department of Mathematics, the University of Hong Kong,
Hong Kong SAR, China}
\email{yujt@hkucc.hku.hk}

\thanks
{The research of Vesselin Drensky was partially supported by Grant
MI-1503/2005 of the Bulgarian National Science Fund.}

\thanks{The research of Jie-Tai Yu was partially
supported by a Hong Kong RGC-CERG Grant.}


\begin{abstract} We study $z$-automorphisms of the polynomial algebra $K[x,y,z]$
and the free associative algebra $K\langle x,y,z\rangle$ over a field $K$, i.e.,
automorphisms which fix the variable $z$. We survey some recent results on such automorphisms and
on the corresponding coordinates. For $K\langle x,y,z\rangle$ we include also results
about the structure of the $z$-tame automorphisms and algorithms which recognize
$z$-tame automorphisms and $z$-tame coordinates.
\end{abstract}

\maketitle

\section*{Introduction}

Let $K$ be an arbitrary field of any characteristic. The automorphism group
$\text{Aut }K[x_1,\ldots,x_n]$ of the polynomial algebra in $n$ variables
is well understood only for $n=1$ and $n=2$. The description is trivial
for $n=1$, when the elements $\varphi\in\text{Aut }K[x_1]$ are defined by
$\varphi(x_1)=\alpha x_1+\beta$, where $\alpha\in K^{\ast}=K\backslash 0$ and
$\beta\in K$. The classical result of Jung \cite{J} and van der Kulk \cite{K}
gives that all automorphisms of $K[x_1,x_2]$ are tame.
Writing the automorphisms of $K[x_1,\ldots,x_n]$ as $n$-tuples
of the images of the variables, and using $x,y$ instead of $x_1,x_2$,
this means that $\text{Aut }K[x,y]$
is generated by the affine automorphisms
\[
\psi=(\alpha_{11}x+\alpha_{21}y+\beta_1,\alpha_{12}x+\alpha_{22}y+\beta_2),
\quad \alpha_{ij},\beta_j\in K,
\]
(and $\psi_1=(\alpha_{11}x+\alpha_{21}y,\alpha_{12}x+\alpha_{22}y)$, the
linear part of $\psi$, is invertible) and the triangular
automorphisms
\[
\rho=(\alpha_1x+p_1(y),\alpha_2y+\beta_2), \quad
\alpha_1,\alpha_2\in K^{\ast},p_1(y)\in K[y],\beta_2\in K.
\]
This result was the starting point of research in several
directions, among them the study of the automorphisms of
$K[x_1,\ldots,x_n]$ for $n>2$, of the automorphisms of the
polynomial algebra $A[x_1,x_2]$ over an arbitrary commutative ring
$A$, and the automorphisms of free associative algebras $K\langle
x_1,\ldots,x_n\rangle$. In all these cases the tame automorphisms
are defined by analogy with the case of $K[x,y]$, as compositions
of affine and triangular automorphisms. One studies not only the
automorphisms but also the coordinates, i.e., the images of $x_1$
under the automorphisms of $K[x_1,\ldots,x_n]$.

Nagata \cite{N} constructed the automorphism of $K[x,y,z]$
\[
\nu=(x-2(y^2+xz)y-(y^2+xz)^2z,y+(y^2+xz)z,z)
\]
which fixes $z$. He showed that $\nu$ is nontame, or wild, considered as an automorphism
of $K[z][x,y]$, and conjectured that it is wild
also as an element of $\text{Aut }K[x,y,z]$.
This example motivated the study in detail of the automorphisms of $K[z][x,y]$ by several
reasons. Their form is simpler than the form of an arbitrary automorphism of $K[x,y,z]$ and
one can apply for their study the results on the automorphisms of $K(z)[x,y]$,
the polynomial algebra in two variables $x,y$, with rational in $z$ coefficients.
Also, the automorphism group of $K[z][x,y]$ provides important examples and conjectures for
$\text{Aut }K[x,y,z]$. We shall mention only few facts related with the topic of the
present paper. It is relatively easy to see (and to decide algorithmically) whether
an endomorphism of $K[z][x,y]$ is an automorphism and whether this automorphism is
$z$-tame, or tame as an automorphism of $K[z][x,y]$, but a similar problem
is much harder for coordinates.
When $K$ is a field of characteristic 0, Drensky and Yu \cite{DY2} presented a simple
algorithm which decides whether a polynomial $f(x,y,z)\in K[x,y,z]$ is a $z$-coordinate
and whether this coordinate is $z$-tame. This provided
many new wild automorphisms and wild coordinates of $K[z][x,y]$.
These results in \cite{DY2} are based on a similar algorithm of Shpilrain and Yu \cite{SY1}
which recognizes the coordinates of $K[x,y]$.
The pioneer work of
Shestakov and Umirbaev \cite{SU1, SU2, SU3}
established that the Nagata automorphism is wild.
It also implies that every wild automorphism of $K[z][x,y]$ is wild as an
automorphism of $K[x,y,z]$.
Umirbaev and Yu \cite {UY} showed that
the $z$-wild coordinates in $K[z][x,y]$ are
wild also in $K[x,y,z]$.
In this way, all $z$-wild examples in \cite{DY2} give automatically
wild examples in $K[x,y,z]$.

Going to free algebras,
Czerniakiewicz \cite{Cz} and Makar-Limanov \cite{ML1, ML2}
proved that all automorphisms of $K\langle x,y\rangle$ are also tame.
There are several candidates for wild automorphisms of
free algebras with more than two generators.
One of them is the example of Anick
$(x+(y(xy-yz),y,z+(zy-yz)y)\in \text{Aut }K\langle x,y,z\rangle$,
see the book by Cohn \cite{C2}, p.~343. It fixes
one variable and its abelianization is a tame automorphism of $K[x,y,z]$.
Exchanging the places of $y$ and $z$, we obtain the automorphism
$(x+z(xz-zy),y+(xz-zy)z,z)$ which fixes $z$ (or a $z$-automorphism),
and refer to it as the Anick automorphism.
It is linear in $x$ and $y$, considering $z$ as a ``noncommutative constant''.
Drensky and Yu \cite{DY3} showed that such $z$-automorphisms are $z$-wild
if and only if a suitable invertible $2\times 2$ matrix with entries from $K[z_1,z_2]$
is not a product of elementary matrices. In particular, this gives that the Anick
automorphism is $z$-wild. For better understanding of $z$-wild automorphisms
one has to know more about $z$-tame automorphisms.
The very recent paper by Drensky and Yu \cite{DY7}
describes the structure of the group of $z$-tame automorphisms of
$K\langle x,y,z\rangle$ and gives algorithms which recognize $z$-tame automorphisms and coordinates
of $K\langle x,y,z\rangle$.

When $K$ is a field of characteristic 0,
Umirbaev \cite{U2} developed further his ideas with Shestakov and described
the defining relations of the group of tame automorphisms of $K[x,y,z]$.
As a result, if $\varphi=(f,g,h)\in \text{Aut }K\langle x,y,z\rangle$ has the property
that the endomorphism $\varphi_0=(f_0,g_0,z)$ of $K\langle x,y,z\rangle$,
where $f_0,g_0$ are the linear in $x,y$ components of $f,g$, respectively,
is a $z$-wild automorphism, then $\varphi$ is wild as an automorphism of
$K\langle x,y,z\rangle$. This implies that the Anick automorphism is wild.
Finally, we want to mention the recent results of Drensky and Yu \cite{DY4, DY5}
which allow to show the wildness of a big class of automorphisms and coordinates of
$K\langle x,y,z\rangle$. Many of them cannot be handled with direct application
of the methods of \cite{DY3} and \cite{U2}.

In the present paper we give a survey of some, mostly
recent results on $z$-automorphisms of $K[x,y,z]$ and $K\langle x,y,z\rangle$.
We present also a selection of
facts on the automorphisms of $K[x,y]$ and $K\langle x,y\rangle$
as a preparation to the study of the $z$-automorphisms
of $K[x,y,z]$ and $K\langle x,y,z\rangle$.
(For more details we refer to the books by van den Essen \cite{E2},
Mikhalev, Shpilrain, and Yu \cite{MSY}, and our survey article \cite{DY1}.)
Finally, we provide a list of open problems and conjectures.

\section{A survey on automorphisms of polynomial algebras}

We fix the field $K$ and consider the polynomial
algebra $K[x,y,z]$ in three variables. We denote its automorphisms
as $\varphi=(f,g,h)$, where $f=\varphi(x)$, $g=\varphi(y)$, $h=\varphi(z)$.
The multiplication is from right to left. If
$\varphi,\psi\in\text{Aut }K[x,y,z]$, then in $\varphi\psi$ we first apply
$\psi$ and then $\varphi$. Hence, if $\varphi=(f,g,h)$ and $\psi=(u,v,w)$, then
\[
\varphi\psi=(u(f,g,h),v(f,g,h),w(f,g,h)).
\]
We call the automorphism $\varphi$ a $z$-automorphism if $\varphi(z)=z$,
and denote
the automorphism group of the $K[z]$-algebra $K[z][x,y]$
by $\text{Aut }K[z][x,y]=\text{Aut}_zK[x,y,z]$. If we want to emphasize
that we work with $z$-automorphisms, we write $\varphi=(f,g)$, omitting the third
coordinate $z$. The affine and triangular automorphisms of $K[x,y,z]$
are, respectively, of the form
\[
\psi=(\alpha_{11}x+\alpha_{21}y+\alpha_{31}z+\beta_1,
\alpha_{12}x+\alpha_{22}y+\alpha_{32}z+\beta_2,
\alpha_{13}x+\alpha_{23}y+\alpha_{33}z+\beta_3),
\]
$\alpha_{ij},\beta_j\in K$, the $3\times 3$ matrix $(\alpha_{ij})$ being invertible,
\[
\rho=(\alpha_1x+p_1(y,z),\alpha_2y+p_2(z),\alpha_3z+\beta_3),
\]
$\alpha_j\in K^{\ast}$, $p_1\in K[y,z]$, $p_2\in K[z]$, $\beta_3\in K$.

We denote by $\text{TAut }K[x,y,z]$ the group of the tame automorphisms which is
generated by affine and triangular automorphisms. When we consider $z$-automorphisms,
the $z$-affine and $z$-triangular automorphisms are, respectively, of the form
\[
\psi=(\alpha_{11}x+\alpha_{21}y+b_1(z),
\alpha_{12}x+\alpha_{22}y+b_2(z)),
\]
again $\alpha_{ij}\in K$, $b_j(z)\in K[z]$,
the $2\times 2$ matrix $(\alpha_{ij})$ being invertible,
\[
\rho=(\alpha_1x+p_1(y,z),\alpha_2y+p_2(z)),
\]
$\alpha_j\in K^{\ast}$, $p_1\in K[y,z]$, $p_2\in K[z]$.
The group of the $z$-tame automorphisms is
$\text{TAut }K[z][x,y]=\text{TAut}_zK[x,y,z]$.
We may define the $z$-affine automorphisms by
\[
\psi=(a_{11}x+a_{21}y+b_1,a_{12}x+a_{22}y+b_2),
\]
where $a_{ij}(z),b_j(z)\in K[z]$ and the $2\times 2$ matrix $(a_{ij})$
is invertible over $K[z]$ (hence its determinant is a nonzero constant in $K$).
But we shall see that this definition is not convenient in the noncommutative
case. For example, the Anick automorphism is affine in this sense but is wild.

We start our survey with the case of two variables accepting notation similar to
the case of three variables.
Jung \cite{J} in 1942 for $K=\mathbb C$
and van der Kulk \cite{K} in 1953 for an arbitrary field $K$ proved that
$\text{Aut }K[x,y]=\text{TAut }K[x,y]$. The proof of van der Kulk gives
also the structure of $\text{Aut }K[x,y]$.

\begin{theorem}\label{Jung-van der Kulk}
Every automorphism of $K[x,y]$ is tame. The group $\text{\rm Aut }K[x,y]$
is isomorphic to the free product $A\ast_CB$ of the group $A$
of the affine automorphisms and the group $B$ of triangular automorphisms
amalgamating their intersection $C=A\cap B$.
\end{theorem}

Every automorphism
$\varphi$ can be presented as a product
\[
\varphi=\psi_m^{\varepsilon_m}\rho_m\psi_{m-1}\cdots\rho_2\psi_1\rho_1^{\varepsilon_1},
\]
where $\psi_i\in A$, $\rho_i\in B$ ($\varepsilon_1$ and $\varepsilon_m$ are
equal to 0 or 1), and, if $\varphi$ does not belong to the union of $A$ and $B$,
we may assume that $\psi_i\in A\backslash B$, $\rho_i\in B\backslash A$.
The freedom of the product means that if $\varphi$ has a nontrivial presentation
of this form, then it is different from the identity automorphism.
Fixing the linear nontriangular automorphism $\tau=(y,x)$,
we can present $\varphi$ in the canonical form
\begin{equation}\label{simplified canonical form of automorphisms}
\varphi=\rho_n\tau\cdots\tau\rho_1\tau\rho_0,
\end{equation}
where $\rho_0,\rho_1,\ldots,\rho_n\in B$ and only $\rho_0$ and $\rho_n$
are allowed to belong to $A$, see for example p. 350 in \cite{C2}.
Using the equalities for compositions of automorphisms
\[
(\alpha x+p(y),\beta y+\gamma)=
(x+\alpha^{-1}(p(y)-p(0)),y)(\alpha x+p(0),\beta y+\gamma),
\quad \gamma\in K,
\]
\[
(\alpha x+\xi,\beta y+\eta)\tau=(\beta y+\eta,\alpha x+\xi)
=\tau(\beta x+\eta,\alpha y+\xi),\quad \xi,\eta\in K,
\]
we can do further simplifications in (\ref{simplified canonical form of automorphisms}),
assuming that $\rho_i=(x+p_i(y),y)$ with $p_i(0)=0$ for all $i=1,\ldots,n$.
We also assume that $\rho_0=(\alpha_0x+p_0(y),\beta_0y+\gamma_0)$.

Let $\varphi=(f,g)$ and let
us assume for example that $\text{deg }p_i(y)=d_i>1$ for all $i=0,1,\ldots,n$,
and let $p_i(y)=\xi_iy^{d_i}+\cdots$, where $\xi_i\in K^{\ast}$
and $q_i(y)=\overline{p_i}(y)=\xi_iy^{d_i}$ is the leading monomial of $p_i(y)$.
Direct computations show that
$\text{deg }f=d_0\cdot \text{deg }g$
and the monomials of highest degree of $f,g$ are, respectively,
\[
\overline{f}=q_0(q_1(\ldots (q_n(y))\ldots))=\xi_0(\zeta y^N)^{d_0},
\]
\[
\overline{g}=\beta_0q_1(\ldots (q_n(y))\ldots)=\beta_0\zeta y^N.
\]
Hence, up to a multiplicative constant, $\overline{f}$ is equal to a power of
$\overline{g}$ and the monomial $q_0(y)$ can be determined uniquely from
$\overline{f},\overline{g}$.
Similar relations hold in the cases when some of the
triangular automorphisms $\rho_0$ and $\rho_n$ is affine but
$\overline{f},\overline{g}$ may be powers of a linear combination of $x,y$
instead of powers of $y$.
These considerations provide an easy algorithm to decide whether an endomorphism of
$K[x,y]$ is an automorphism, see Theorem 6.8.5 in \cite{C2}. Of course, in all
algorithms we assume that the field $K$ is constructive.

\begin{algorithm}\label{algorithm for recognizing of automorphisms}
Let $\varphi=(f,g)$ be an endomorphism of $K[x,y]$.

{\it Step 0}. If some of the polynomials $f,g$ is a constant from $K$, then
$\varphi$ is not an automorphism.

{\it Step 1}. Let $u,v$ be the homogeneous components of highest degree of $f,g$,
respectively.
If both $u,v$ are linear, then we check whether they are linearly
independent. If yes, then $\varphi$ is an affine automorphism. If
$u,v$ are linearly dependent, then $\varphi$ is not an automorphism.

{\it Step 2}. If $u$ is not linear,
$\text{deg }u\geq \text{deg }v$, and $u=\beta v^d$ for some
$\beta\in K^{\ast}$ and some $d\geq 1$, then we replace
$\varphi=(f,g)$ with $\varphi_1=(f-\beta g^d,g)$.
Then we apply Step 0 to $\varphi_1$.
If $u$ cannot be presented in the form $u=\beta v^d$,
then $\varphi$ is not an automorphism.

{\it Step 3}. If $v$ is not linear,
$\text{deg }u<\text{deg }v$, and $v=\beta u^d$, then we replace
$\varphi=(f,g)$ with $\varphi_1=(f,g-\beta f^d)$.
Then we apply Step 0 to $\varphi_1$.
If $v$ cannot be presented in the form $v=\beta u^d$,
then $\varphi$ is not an automorphism.

Steps 2 and 3 of the algorithm correspond to a triangular
automorphism and the algorithm also provides a decomposition of
$\varphi$ as a product of affine and triangular automorphisms.
If $\varphi_1=(f-\beta g^d,g)$ in Step 2, and
$\sigma_1=(x-\beta y^d,y)$, then
$\varphi_1=\varphi\sigma_1$ and $\varphi=\varphi_1\sigma_1^{-1}$.
If $\varphi_1=(f,g-\beta f^d)$
in Step 3, then $\varphi_1=\varphi\tau\sigma_1\tau$ and
$\varphi=\varphi_1\tau\sigma_1^{-1}\tau$.
\end{algorithm}

In the general case of polynomial algebras in several variables, there is
an effective algorithm which decides whether an endomorphism is an automorphism.
It involves Gr\"obner bases techniques, see our survey article \cite{DY1}
for references.

Algorithm \ref{algorithm for recognizing of automorphisms} can be modified to decide
whether a polynomial $f(x,y)$ is a coordinate of $K[x,y]$.
We start with the analysis of the behavior of the first coordinate $f$
of $\varphi$ in (\ref{simplified canonical form of automorphisms}).
Note, that if $\varphi=(f,g)$ and $\varphi'=(f,g')$ are two automorphisms with the
same first coordinate, then $\varphi^{-1}\varphi'$ fixes $x$. Hence
$\varphi^{-1}\varphi'=(x,g'')$ and, obligatorily, $g''=\beta y+r(x)$.
In this way, if we find one ``coordinate mate'' $g$ of $f$, then we can find
all other mates.
Let $(a,b)$ be a pair of positive integers. We define the
$(a,b)$-degree of a monomial $x^my^n$ as $\text{deg}_{(a,b)}x^my^n=am+bn$,
and denote the homogeneous component of maximal
$(a,b)$-degree of $f\in K[x,y]$ by $\vert f\vert_{(a,b)}$.
Let us assume again that $\text{deg }p_i(y)>1$ for all $i=0,1,\ldots,n$, and
let $h$ be the first coordinate of
$\psi=\rho_{n-1}\tau\cdots\tau\rho_1\tau\rho_0$. Then
\[
\overline{h}=q_0(q_1(\ldots (q_{n-1}(y))\ldots))=\vartheta y^M,\quad
\vartheta\in K^{\ast}, M\geq 2.
\]
The homogeneous component of maximal $(d_n,1)$-degree of $x+p_n(y)$
is $x+\xi_ny^{d_n}$ and direct calculations give
\[
\vert f\vert_{(d_n,1)}=\vert \rho_n\tau(\overline{h})\vert_{(d_n,1)}
=\vartheta (x+\xi_ny^{d_n})^M.
\]
Hence we can recover from here $d_n$ and $\xi_n$.
In particular,
\[
f=\omega(x^M+\xi_n^My^{Md_n})+\cdots.
\]
The considerations are similar when
some of the automorphisms $\rho_0$ and $\rho_n$ is affine.
We easily obtain the validity of the following algorithm which recognizes
the coordinates in $K[x,y]$.

\begin{algorithm}\label{algorithm for recognizing of coordinates}
Let $f(x,y)\in K[x,y]$ be a polynomial.

{\it Step 0}. If $f$ is a constant from $K$, then
$f$ is not a coordinate. If $f=\lambda_1 x+\lambda_2 y+\lambda_3$ is linear, then
$f$ is a coordinate. One of its mates is of the form
$g=\mu_1 x+\mu_2 y+\mu_3$ and can be determined from the property
that the polynomials $\lambda_1 x+\lambda_2 y$ and $\mu_1 x+\mu_2 y$
are linearly independent.

{\it Step 1}.
If $\text{deg }f>1$ and $f$ depends only on $x$ (or only on $y$),
then it is not a coordinate. If $f$ does not contain as a summand
$\eta x^M$ or $\zeta y^N$ for some $\eta,\zeta\in K^{\ast}$ and $M,N\geq 1$,
then $f$ is not a coordinate.

{\it Step 2}. Let $f=\eta x^M+\zeta y^N+\cdots$,
where $\cdots$ stays for the linear combination of the monomials of the form
$x^j$, $j<M$, $y^k$, $k<N$, and the monomials $x^ry^s$, $r,s\geq 1$.
If $M$ does not divide $N$ and $N$ does not divide $M$, then $f$ is not
a coordinate.

{\it Step 3}. Let $M$ divide $N$ and $N=Md$. Write $\zeta$ in the form
$\zeta=\eta\xi^d$, $\xi\in K^{\ast}$. (If $K$ is not algebraically closed
and this is not possible, then $f$ is not a coordinate of $K[x,y]$.)
Consider the $(d,1)$-grading of $K[x,y]$.
If $\vert f\vert_{(d,1)}\not=\eta(x+\xi y^d)^M$, then $f$ is not a coordinate.
If $\vert f\vert_{(d,1)}=\eta(x+\xi y^d)^M$ and $M>1$, then
define $\vartheta=(x+\xi y^d,y)$, replace
$f$ with $f_1=\vartheta^{-1}(f)$
and go to Step 0.
If $M=1$ and $\vert f\vert_{(d,1)}=\eta(x+\xi y^d)$, then
$u=f-\eta(x+\xi y^d)$ does not depend on $x$ and hence
$f=\eta x+\xi y^d+u(y)$ is a coordinate. As a mate, we can take $g=y$.

{\it Step 4}. Let $N$ divide $M$ and $M=Nd$, $d>1$. The considerations are
similar to those of Step 3, working with $\vartheta=(x,y+\xi x^d)$
instead of with $\vartheta=(x+\xi y^d,y)$.

As in the previous algorithm,
Steps 3 and 4 correspond to a triangular
automorphism. For example, if
$\vert f\vert_{(d,1)}=\eta(x+\xi y^d)^M$ and $M>1$ in Step 3,
and $u=f-\eta(x+\xi y^d)^M$, then
$\vartheta^{-1}(f)=\vartheta^{-1}(u)+\eta x^M$. Since
$x+\xi y^d$ and $y$ are $(d,1)$-homogeneous, the $(d,1)$-degree of
$\vartheta^{-1}(u)$ is smaller than this of $f$. Hence
$\vartheta^{-1}(f)$ does not contain summands of the form $\zeta_jy^j$,
$j\geq Md=N$ and we decrease the sum $M+N$.
\end{algorithm}

To the best of our knowledge, the above algorithm had not been explicitly stated
before \cite{SY3}, where Shpilrain and Yu established an algorithm
which gives a canonical form, up to automorphic equivalence,
of a class of polynomials in $K[x,y]$. The automorphic equivalence problem
for $K[x,y]$ asks how to decide whether, for two given polynomials $p,q\in K[x,y]$,
there exists an automorphism $\varphi$ such that $q=\varphi(p)$. It was solved
over $\mathbb C$ by Wightwick \cite{Wi} and, over an arbitrary algebraically closed
constructive field $K$, by Makar-Limanov, Shpilrain, and Yu \cite{MLSY}.

For a polynomial $f(x,y)\in K[x,y]$ we denote by
\[
f_x=\frac{\partial f}{\partial x},\quad f_y=\frac{\partial f}{\partial y}
\]
the partial derivatives of $f$ with respect to $x$ and $y$.
If $\varphi=(f,g)$ is an endomorphism of $K[x,y]$, then its Jacobian matrix is
\[
J(\varphi)=\left\vert
\begin{matrix}
f_x&g_x\\
f_y&g_y\\
\end{matrix}
\right\vert.
\]
For polynomial algebras in more than two variables the Jacobian matrix
is defined in a similar way.
The chain rule gives that
\[
J(\varphi\psi)=J(\varphi)\varphi\left(J(\psi)\right),
\]
where $\varphi$ acts on $J(\psi)$ componentwise. This implies that if $\varphi$ is an
automorphism, then its Jacobian matrix $J(\varphi)$
is invertible and the determinant $\text{det}J(\varphi)$ is a nonzero scalar in $K$.
The famous Jacobian conjecture states that if $\text{char }K=0$ and $J(\varphi)$
is invertible, then the endomorphism $\varphi$ is an automorphism.

If $\rho$ is a triangular automorphism, then its Jacobian matrix
$J(\rho)$ is a low triangular matrix. If $\psi$ is an affine automorphism, then
$J(\psi)$ belongs to $GL_2(K)$ and also is a product of low and upper triangular matrices.
Since the automorphisms of $K[x,y]$ are tame,
the chain rule gives that $J(\varphi)$ is a product of triangular matrices
for all $\varphi\in\text{Aut }K[x,y]$. Hence $J(\varphi)$ belongs to the group
$GE_2(K[x,y])$ generated by the elementary $2\times 2$ matrices with entries
from $K[x,y]$ and diagonal matrices. This fact was used by Wright \cite{Wr}
to solve the Jacobian conjecture for endomorphisms $\varphi$ of $K[x,y]$ with
the property $J(\varphi)\in GE_2(K[x,y])$. Wright \cite{Wr} established also
that the group $GE_2(K[x_1,\ldots,x_n])$ is a free product
of $GL_2(K)$ and the group of low triangular matrices with entries from
$K[x_1,\ldots,x_n]$ with amalgamation of their intersection.
Hence $GE_2(K[x_1,\ldots,x_n])$ has a similar structure as $\text{Aut }K[x,y]$.
There is an algorithm which decides whether a matrix
in $GL_2(K[x_1,\ldots,x_n])$ belongs to $GE_2(K[x_1,\ldots,x_n])$.
It was suggested by Tolhuizen, Hollmann and Kalker \cite{THK}
for the partial ordering by degree and then generalized  by Park \cite{P}
for any term-ordering on $K[x_1,\ldots,x_n]$. One applies Gaussian
elimination process on the matrix based on the Euclidean division algorithm
for $K[x_1,\ldots,x_n]$. The matrix belongs to $GE_2(K[x_1,\ldots,x_n])$
if and only if this procedure brings it to an elementary or diagonal matrix.
The idea of Wright \cite{Wr} and this algorithm were combined by
Shpilrain and Yu \cite{SY1} to give a very simple algorithm
which decides whether a polynomial
$f(x,y)\in K[x,y]$, $\text{char }K=0$, is a coordinate.

Let $K$ be a constructive field of characteristic 0.
We fix a term-ordering on $K[x,y]$. Recall that this
a linear order on the set of monomials $x^my^n$,
which is a multiplicatively compatible well-ordering.

\begin{algorithm}\label{algorithm for coordinates in two variables}
(Shpilrain and Yu \cite{SY1})
Let $f(x,y)\in K[x,y]$.
Consider the partial derivatives $p=f_x$ and $q=f_y$.

{\it Step 0}. If $p=q=0$, then $f$ is not a coordinate.
If either $p=0$ and $q\in K^{\ast}$, or $q=0$ and $p\in K^{\ast}$
then $f$ is a coordinate.
If either $p=0$ and $\text{deg }q>0$, or $q=0$ and $\text{deg }p>0$
then $f$ is not a coordinate.

{\it Step 1}. If $p,q\not=0$, apply the
Euclidean division algorithm to $p$ and $q$. If this is impossible, then $f$ is not
a coordinate. If $p=qs+r$, where $r,s\in K[x,y]$ with $r=0$ or
$\text{deg }r<\text{deg }q$, then replace $p$ with $r$ and go to Step 0.
If $q=ps+r$, where $r,s\in K[x,y]$ with $r=0$ or
$\text{deg }r<\text{deg }p$, then replace $q$ with $r$ and again go to Step 0.
\end{algorithm}

Before \cite{SY1, SY3}, more complicated algorithms which decide whether
$f(x,y)\in{\mathbb C}[x,y]$ is a coordinate were given by
Ch\c adzy\'nski and Krasi\'nski \cite{CK} in terms of the so called
\L ojasiewicz exponent at infinity
and by van den Essen \cite{E1} using the technique of locally nilpotent derivations.

Now we continue our survey with exposition on automorphisms of $K[z][x,y]$ and $K[x,y,z]$.

The first example of a $z$-wild automorphism of $K[z][x,y]$
was constructed by Nagata \cite{N}.

\begin{theorem}\label{Nagata automorphism is z-wild}
{\rm (Theorem 1.4, p. 42 \cite{N})}
The Nagata automorphism
\[
\nu=(x-2y(y^2+xz)-z(y^2+xz)^2,y+z(y^2+xz))
\]
of $K[z][x,y]$ is $z$-wild.
\end{theorem}

\begin{conjecture}\label{Nagata conjecture}
The Nagata automorphism is wild as an automorphism of
the polynomial algebra $K[x,y,z]$.
\end{conjecture}

The Nagata automorphism is among the most investigated automorphisms of
$K[x,y,z]$. It was used as a counterexample to different problems
on automorphisms of polynomial algebras and group actions on affine
spaces, see e.g. the paper by Bass \cite{Ba} and our survey \cite{DY1}.
We shall consider only one property of the Nagata automorphism.

Recall that a $K$-linear operator $\delta$ acting on $K[x_1,\ldots,x_n]$
is a derivation, if $\delta(uv)=\delta(u)v+u\delta(v)$
for all $u,v\in K[x_1,\ldots,x_n]$. The derivation $\delta$
is triangular if $\delta(x_j)\in K[x_{j+1},\ldots,x_n]$, $j=1,\ldots,n$.
The derivation $\delta$ is locally nilpotent if for every $u$ there exists a $d$
such that $\delta^d(u)=0$. If $\delta$ is locally nilpotent
and $\text{char }K=0$, then the formal series
\[
\text{exp}(\delta)=1+\frac{\delta}{1!}+\frac{\delta^2}{2!}+\cdots
\]
is a well defined linear operator on $K[x_1,\ldots,x_n]$ which is an
automorphism of the polynomial algebra.
Clearly, all triangular derivations are locally nilpotent.
The corresponding exponential automorphisms are also triangular and hence tame.
Another way to construct locally nilpotent derivations is to start with
a locally nilpotent derivation $\delta$ and to consider its kernel $\text{Ker }\delta$.
If $0\not=w\in \text{Ker }\delta$, then $\Delta=w\delta$ is also a
locally nilpotent derivation with the same kernel.

\begin{example}\label{Nagata automorphism using derivations}
Let $\text{char }K=0$ and let $\delta$ be the triangular derivation of $K[x,y,z]$
defined by
\[
\delta(x)=-2y,\quad \delta(y)=z,\quad \delta(z)=0.
\]
Then $\text{Ker }\delta=K[y^2+xz,z]$. For $w=y^2+xz$ and $\Delta=w\delta$
we obtain that $\nu=\text{exp}(\Delta)$, the Nagata automorphism, because
\[
\text{exp}(\Delta):x\to (1+\frac{\Delta}{1!}+\frac{\Delta^2}{2!}+\cdots)(x)
\]
\[
= x+(y^2+xz)\frac{\delta(x)}{1!}+(y^2+xz)^2\frac{\delta^2(x)}{2!}+\cdots
=x-2(y^2+xz)y-(y^2+xz)^2z,
\]
\[
\text{exp}(\Delta):y\to (1+\frac{\Delta}{1!}+\frac{\Delta^2}{2!}+\cdots)(y)
\]
\[
= y+(y^2+xz)\frac{\delta(y)}{1!}+\cdots
=y+(y^2+xz)z,
\]
\[
\text{exp}(\Delta):z\to (1+\frac{\Delta}{1!}+\frac{\Delta^2}{2!}+\cdots)(z)=z.
\]
\end{example}

The automorphism $\varphi$ of the polynomial algebra $K[x_1,\ldots,x_n]$
is called stably tame, if it becomes tame when extended identically, by
$\varphi(x_{n+i})=x_{n+i}$, $i=1,\ldots,m$,
on the polynomial algebra $K[x_1,\ldots,x_{n+m}]$ for a suitable integer $m$.
The following theorem of Smith \cite{Sm} gives a big class of stably tame automorphisms.

\begin{theorem}\label{theorem of Martha Smith}
Over a field $K$ of characteristic $0$,
every triangular derivation
$\delta$ of $K[x_1,\ldots,x_n]$ and every $w\in\text{\rm Ker }\delta$,
the automorphism $\exp(w\delta)$ is stably tame and becomes
tame on $K[x_1,\ldots,x_n,x_{n+1}]$ assuming that $\delta(x_{n+1})=0$.
\end{theorem}

The key observation in the proof of the theorem is the
following equation in $K[x_1,\ldots,x_n,x_{n+1}]$ discovered by Smith
\[
\text{exp}(w\delta)=
\vartheta^{-1}\cdot\text{exp}(x_{n+1}\delta)\cdot\vartheta\cdot\text{exp}^{-1}(x_{n+1}\delta),
\]
where $\vartheta=(x_1,\ldots,x_n,x_{n+1}+w)$.
Since $w\in K[x_1,\ldots,x_n]$, the automorphism $\vartheta$
is tame (triangular assuming that $x_{n+1}<x_i$, $i=1,\ldots,n$),
the automorphism $\text{exp}(x_{n+1}\delta)$ is also
tame because $\delta$ is a triangular derivation
(and $\exp(x_{n+1}\delta)$ is a triangular automorphism assuming that
$x_i<x_{n+1}$). Hence, $\exp(w\delta)$ is also tame as a composition
of tame automorphisms of $K[x_1,\ldots,x_n,x_{n+1}]$.

The approach of Smith is one of the main steps in the proofs of all known results on
stable tameness of automorphisms of polynomial algebras. Recently, the notion
of stable tameness was transferred also to coordinates, see Berson \cite{B1, B2}.
The coordinate $f\in K[x_1,\ldots,x_n]$ is called stably tame, if,
for some $m$, there exists a tame automorphism $\varphi$ of $K[x_1,\ldots,x_{n+m}]$
such that $\varphi(x_1)=f$. This tameness is weaker than the tameness of
automorphisms because one does not require that $\varphi$ fixes the variables
$x_{n+1},\ldots,x_{n+m}$. Berson \cite{B2} showed that, over a noetherian
$\mathbb Q$-domain $A$ of Krull dimension one, all coordinates of $A[x,y]$
are stably tame. Hence every $z$-coordinate of $K[z][x,y]$ is stably tame.
For a class of coordinates of $A[x,y]$ Edo \cite{Ed} proved that they are
coordinates of stably tame automorphisms. We shall mention the following
well known problem.

\begin{problem}
Are all automorphisms of $K[z][x,y]$ stably tame?
\end{problem}

Every $z$-automorphism of $K[z][x,y]$ is also an automorphism of the polynomial algebra
$K(z)[x,y]$ over the field of fractions $K(z)$. Hence one may use the theory of
automorphisms over a field to study $z$-automorphisms.
Let us consider the chain of inclusions
\[
\text{TAut }K[z][x,y]\subset\text{Aut }K[z][x,y]
\subset \text{Aut }K(z)[x,y]
\]
and the equality
\[
\text{Aut }K(z)[x,y]=A(K(z)[x,y])\ast_{C(K(z)[x,y])}B(K(z)[x,y]),
\]
where $A(K(z)[x,y])$, $B(K(z)[x,y])$, and $C(K(z)[x,y])$ are, respectively,
the groups of affine and triangular automorphisms of $K(z)[x,y]$, and
the intersection of these two groups. The nice structure of
$\text{Aut }K(z)[x,y]$ should imply some structure on $\text{TAut }K[z][x,y]$
and on $\text{Aut }K[z][x,y]$. For any commutative ring
$A$ we denote by $\text{Aut}^0A[x,y]$ the subgroup of $\text{Aut }A[x,y]$
consisting of all automorphisms $\varphi=(f,g)$ with the property
$f(0,0)=g(0,0)=0$. We use similar notation for other groups of automorphisms.
We call the elements of
$\text{Aut}^0A[x,y]$ automorphisms without constant terms.
For many purposes the study of $\text{Aut}^0A[x,y]$ is equivalent
to that of $\text{Aut }A[x,y]$ because every automorphism
of $A[x,y]$ is a composition of a translation $(x+\alpha,y+\beta)$ and
an automorphism without constant terms. Although the following result of Wright \cite{Wr}
holds in a more general situation, we shall state it in the case of $K[z][x,y]$ only.

\begin{theorem}\label{theorem of Wright}
{\rm (i)}
Over any field $K$ the group $\text{\rm TAut }K[z][x,y]$ and its subgroup
$\text{\rm TAut}^0K[z][x,y]$ have the amalgamated free product structure
\[
\text{\rm TAut }K[z][x,y]=A\ast_CB,\quad
\text{\rm TAut}^0K[z][x,y]=GL_2(K)\ast_{C^0}B^0,
\]
where $A$ and $B$ are, respectively, the group of affine $z$-automorphisms
and the group of triangular $z$-automorphisms and $C=A\cap B$.

{\rm (ii)}
There is a subgroup $W$ of $\text{\rm Aut}^0K[z][x,y]$ containing
the group $B^0(K[z][x,y])$ of all triangular automorphisms
without constant terms such that
\[
GL_2(K[z])\cap W=GL_2(K[z])\cap B^0(K[z][x,y])=C^0(K[z][x,y]),
\]
and
\[
\text{\rm Aut}^0K[z][x,y]\cong GL_2(K[z])\ast_{C^0(K[z][x,y])}W,
\]
the amalgamated free product of $GL_2(K[z])$ and $W$.
\end{theorem}

The original statement of Theorem \ref{theorem of Wright} involves affine
and linear automorphisms with coefficients from $K[z]$ but this is not
essential because every invertible matrix with entries in $K[z]$ is a product
of elementary and diagonal matrices.

Theorem \ref{theorem of Wright} allows to obtain immediately analogues of
Algorithms \ref{algorithm for recognizing of automorphisms}
and \ref{algorithm for recognizing of coordinates} which recognize, respectively,
the $z$-tame automorphisms and coordinates, and find the decomposition
of the automorphism as a product of affine and triangular automorphisms.
Algorithm \ref{algorithm for coordinates in two variables}
also has its analogue for $z$-tame coordinates of $K[z][x,y]$.

\begin{proposition}\label{algorithm for z-tame coordinates}
{\rm (Drensky and Yu \cite{DY2})}
Let $K$ be a constructive field of characteristic $0$ and let
$f(x,y,z)\in K[x,y,z]$.
Consider the partial derivatives $p=f_x$ and $q=f_y$. Then the procedure
described in Algorithm \ref{algorithm for coordinates in two variables}
decides whether $f$ is a $z$-tame coordinate of $K[z][x,y]$.
\end{proposition}

Of course, not every coordinate of $K(z)[x,y]$ which is in $K[z][x,y]$
is a coordinate also in $K[z][x,y]$.
Now we shall describe an algorithm which handles this problem.
Recall that $f\in K[z][x,y]$ has a unimodular gradient
with respect to $x$ and $y$ if the partial derivatives
$f_x$ and $f_y$ generate $K[z][x,y]$ as an ideal.
The following result is obtained by Drensky and Yu \cite{DY2}.
Their proof uses a result of
Daigle and Freudenburg \cite{DF} on locally nilpotent derivations.
The same proof works in the general case, over an arbitrary
commutative $\mathbb Q$-algebra $A$, see Edo and V\'en\'ereau \cite{EdV}.

\begin{theorem}\label{description of any z-coordinate}
{\rm (Drensky and Yu \cite{DY2})}
The polynomial $f(x,y)\in K[z][x,y]$ is a coordinate
in $K[z][x,y]$ if and only if it is with unimodular gradient in $K[z][x,y]$
and is a coordinate in $K(z)[x,y]$.
\end{theorem}

Combined with Algorithm \ref{algorithm for coordinates in two variables},
Theorem \ref{description of any z-coordinate}
allows to solve
effectively the problem whether a polynomial is a coordinate, but this
time we have to apply Gr\"obner bases techniques instead of the Euclidean
algorithm.

\begin{algorithm}\label{algorithm for z-coordinates}
{\rm (Drensky and Yu \cite{DY2})}
Given a polynomial $f\in K[z][x,y]$, where $K$ is a constructive field of
characteristic $0$,
we want to decide whether it is
a coordinate.

{\it Step 0}. Take the partial derivatives
$f_x, f_y\in K[z][x,y]$.

{\it Step 1}. Fix some term-ordering in $K[x,y,z]$.
Find the reduced Gr\"obner basis of the ideal
of $K[x,y,z]$ generated by the polynomials $f_x,f_y$.
This ideal coincides with $K[x,y,z]$
if and only if the obtained reduced Gr\"obner basis consists
of a nonzero constant in $K$ only. Hence,
if the Gr\"obner basis does not consist of
a nonzero constant, then $f$
is not a coordinate polynomial of $K[z][x,y]$.

{\it Step 2}. Working on $K(z)$ instead on $K[z]$
(and fixing a term-ordering in $K[z][x,y]$), we apply
Algorithm \ref{algorithm for coordinates in two variables}
and determine whether $f$ is a coordinate in
$K(z)[x,y]$. If the answer is negative,
then $f$ cannot be a coordinate in $K[z][x,y]$.
Otherwise $f$ is a coordinate of $K[z][x,y]$.
\end{algorithm}

For example, if $f=y+(y^2+xz)z$ is the second coordinate of the Nagata automorphism,
then $f_x=z^2$, $f_y=1+2yz$. It is easy to see that $f_x$ and $f_y$ generate
as an ideal the whole $K[x,y,z]$. (For the proof, $1=-4y^2f_x+(1-2yz)f_y$
belongs to the ideal generated by $f_x$ and $f_y$.)
Applying the Euclidean algorithm in $K(z)[x,y]$,
we see that $f_y$ is divisible by $f_x$, and $f_x$ is a constant,
as an element of $K(z)$. Hence $f$ is a coordinate both in $K(z)[x,y]$ and $K[z][x,y]$.
Since the leading monomials $z^2$ and $2yz$ of $f_x$ and $f_y$
are not divisible by each other, we obtain that $f$ is a $z$-wild coordinate.

See also our survey \cite{DY1} and
the paper by Berson and van den Essen \cite{BE}
for the case of coordinates in $A[x,y]$ over any finitely generated
$\mathbb Q$-algebra $A$.
Going back to the automorphisms of $K[x,y,z]$,
Theorem \ref{description of any z-coordinate}
and Algorithm \ref{algorithm for z-coordinates}
allow to find effectively
a lot of wild automorphisms
of $K[z][x,y]$, giving in this way new candidates
for wild automorphisms of
$K[x,y,z]$, all of them fixing $z$ as
in the example suggested by Nagata, see \cite{DY1, DY2}.
One of the possible ways to search for new $z$-wild automorphisms is
to try to find triangular automorphisms $\rho_0,\rho_1,\ldots,\rho_n$ of $K(z)[x,y]$
such that $\rho_0$ and $\rho_n$ are not automorphisms of $K[z][x,y]$
but, nevertheless, $\varphi=\rho_n\tau\cdots\tau\rho_1\tau\rho_0$
in the canonical form (\ref{simplified canonical form of automorphisms})
belongs to $\text{Aut }K[z][x,y]$ (and not to $\text{TAut }K[z][x,y]$).
Such a search was performed for small $n$. For example, the Nagata automorphism,
considered as an automorphism of $K(z)[x,y]$ has the presentation
\begin{equation}\label{Nagata automorphism in conjugate form}
\nu=\rho_0\rho_1\rho_0^{-1},\quad\text{where}\quad
\rho_0=\left(x+\frac{y^2}{z},y\right),\quad
\rho_1=(x,y+z^2x).
\end{equation}

There are also other ways to construct automorphisms of polynomial algebras.
One of them uses the algebra $R_{22}$ generated by two generic $2\times 2$ matrices
\[
X=\left(
\begin{matrix}
x_{11}&x_{12}\\
x_{21}&x_{22}\\
\end{matrix}
\right),\quad
Y=\left(
\begin{matrix}
y_{11}&y_{12}\\
y_{21}&y_{22}\\
\end{matrix}
\right),
\]
where $x_{ij},y_{ij}$ are commuting variables. Let $T_{22}$ be
the noncommutative trace algebra generated by $R_{22}$ and
all traces $\text{tr}(u)$, $u\in R_{22}$.
If $\text{char }K\not=2$, then the center $C(T_{22})$ of $T_{22}$ is isomorphic to
the polynomial algebra generated by the
five algebraically independent polynomials
\[
\text{tr}(X),\text{tr}(Y),\text{tr}(X^2),\text{tr}(XY),\text{tr}(Y^2).
\]
Hence every automorphism of the generic matrix algebra $R_{22}$ induces
automorphisms on $T_{22}$ and $C(T_{22})\cong K[x_1,\ldots,x_5]$.
This construction was used to obtain different automorphisms of $K[x_1,\ldots,x_5]$,
see Bergman \cite{Be}, Alev and Le Bruyn \cite{AL}, Drensky and Gupta \cite{DG1, DG2}.
Nevertheless, it has turned out that all these automorphisms can be described
in terms of locally nilpotent derivations, and as naturally looking
automorphisms of $K(y_1,\ldots,y_m)[x_1,\ldots,x_n]$.
Compare, for example, the results of \cite{DG1} with those
of Drensky, van den Essen, and Stefanov \cite{DES}.

Till the end of the section we assume that the field $K$ is of characteristic 0.
Shestakov and Umirbaev \cite{SU1, SU2, SU3} developed a special technique
and established, more than 30 years after Nagata discovered his automorphism,
that the Nagata automorphism is wild, as an automorphism of $K[x,y,z]$.
A popular exposition of the ideas of their proof is given by van den Essen \cite{E3}.
We shall state some of the results in \cite{SU1, SU2, SU3}
and shall say a couple of words for the (very rough) idea of the proofs.

\begin{theorem}\label{theorem of Shestakov and Umirbaev}
{\rm (Shestakov and Umirbaev \cite{SU1, SU2, SU3})}
Over a constructive field $K$ of characteristic $0$, there exists
an algorithm which decides whether an automorphism of $K[x,y,z]$ is tame.
\end{theorem}

If $\varphi=(f,g)$ is an automorphism of the polynomial algebra in two variables,
which by Theorem \ref{Jung-van der Kulk} is tame,
then Algorithm \ref{algorithm for recognizing of automorphisms}
allows to find a triangular or a linear automorphism which, applied to $f$ and $g$,
decreases the sum of their degrees. Hence, we may find a sequence
$\psi_1,\ldots,\psi_n$ of affine and triangular automorphisms such that,
for $\varphi_i=\psi_i\cdots\psi_1\varphi$ we have
$\text{deg}(\varphi_{i+1}(x))+\text{deg}(\varphi_{i+1}(y))
<\text{deg}(\varphi_i(x))+\text{deg}(\varphi_i(y))$.
But in the case of $K[x,y,z]$ it may happen that this is not possible.
The main idea of the proof of Shestakov and Umirbaev is to use a peak reduction.
If $\varphi=(f,g,h)$ is a tame automorphism of $K[x,y,z]$,
then one follows the sum $\text{deg }f+\text{deg }g+\text{deg }h$ and tries to
minimize it. Let $\overline{f}$ be the homogenous component of maximal degree of $f$.
If one of the polynomials
$\overline{f},\overline{g},\overline{h}$
belongs to the subalgebra of $K[x,y,z]$ generated by the other two,
say $\overline{f}=p(\overline{g},\overline{h})$,
then one can reduce the degree immediately replacing $\varphi$ with the
tame automorphism $\varphi(x-p(y,z),y,z)$. If neither of
the polynomials $\overline{f},\overline{g},\overline{h}$
belongs to the subalgebra generated by the others, then the degree cannot be decreased
immediately. Nevertheless, depending on the type of
the triple $(\overline{f},\overline{g},\overline{h})$
the authors find a sequence of several affine and triangular automorphisms which
decreases the degree. For this purpose, they need a detailed information on
the two-generated subalgebras of $K[x,y,z]$ and lower estimates for the degree
of the elements of these subalgebras. This is achieved by embedding the polynomial
algebra in the free Poisson algebra (or the algebra of universal Poisson brackets)
and the usage systematically the brackets as an additional tool.
So, the proof involves essentially methods of noncommutative (and even nonassociative) algebra.

Applying the algorithm of Theorem \ref{theorem of Shestakov and Umirbaev}
Shestakov and Umirbaev solved into affirmative the Nagata conjecture.

\begin{theorem}\label{Nagata automorphism is wild}
{\rm (Shestakov and Umirbaev \cite{SU1, SU2, SU3})}
The Nagata automorphism is wild as an automorphism of $K[x,y,z]$.
\end{theorem}

They also obtained the following result.

\begin{theorem}\label{z-wild means wild}
Every $z$-automorphism of $K[x,y,z]$ which is $z$-wild
is also wild as an automorphism of $K[x,y,z]$.
\end{theorem}

This theorem replaces the difficult problem
to decide whether a $z$-automorphism is wild as an automorphism of $K[x,y,z]$
with the easier one whether it is $z$-wild. Now, one uses the results of
Drensky and Yu \cite{DY2} to find new examples of wild automorphisms and to determine
whether a $z$-automorphism is wild.

Following Umirbaev and Yu \cite{UY}, a coordinate $f$ in $K[x,y,z]$
is called wild, if every automorphism $\varphi=(f,g,h)$
having $f$ as a first coordinate is wild. Umirbaev and Yu proved the following theorem
which solves the Strong Nagata Conjecture.

\begin{theorem}\label{strong Nagata conjecture}
The algebra $K[x,y,z]$ possesses wild coordinates. If $\varphi=(f,g,z)$
is a $z$-wild automorphism of $K[x,y,z]$, then its nontrivial coordinates $f,g$
are both wild.
\end{theorem}

The direct application of the algorithm of Shestakov and Umirbaev
from Theorem \ref{theorem of Shestakov and Umirbaev} works successfully
for $z$-automorphisms. It is not clear how it may give explicit wild automorphisms
which do not fix $z$. Theorem \ref {strong Nagata conjecture}
of Umirbaev and Yu gives new examples of wild automorphisms.

\begin{example}
Let $\nu$ be the Nagata automorphism and let $\tau_{x,z}=(z,y,x)$.
The automorphism $\varphi=\nu\tau_{x,z}\nu$ does not fix any of the variables
$x,y,z$. Nevertheless it is wild because $\varphi(z)=\nu(x)$
and $\nu(x)$ is a wild coordinate.
\end{example}

It would be interesting to give examples of wild automorphisms
of $K[x,y,z]$ which cannot be obtained from Theorems
\ref{z-wild means wild},
\ref{strong Nagata conjecture}, and \ref{description of any z-coordinate},
applying Algorithm \ref{algorithm for z-coordinates}.
There are several candidates for such wild automorphisms.
All they come from the construction of Freudenburg \cite{F1, F2}
of locally nilpotent derivations
which do not annulate any coordinate of $K[x,y,z]$.
Among them are derivations which annulate the form $y^2+xz$. The simplest
derivation of this kind is
defined as the determinant of a Jacobian matrix
\[
\Delta(u)=\left\vert
\begin{matrix}
f_x&g_x&u_x\\
f_y&g_y&u_y\\
f_z&g_z&u_z\\
\end{matrix}
\right\vert,
\]
where the polynomials $f$ and $g$ are given by
\[
f=y^2+xz,\quad g=zf^2+2x^2yf-x^5.
\]
Since $\Delta$ is locally nilpotent, $\text{exp}(\Delta)$ is an automorphism.
It does not fix any coordinate, so Theorems \ref{z-wild means wild} and
\ref{strong Nagata conjecture} cannot be applied directly.

\begin{problem}
Is the automorphism $\text{\rm exp}(\Delta)$ defined above a wild automorphism of
$K[x,y,z]$?
\end{problem}

We want to mention, that for $n>3$ no wild automorphisms of
$K[x_1,\ldots,x_n]$ are known.

The isomorphism $\text{Aut }K[x,y]\cong A\ast_CB$ means that a set of defining relations of
$\text{Aut }K[x,y]$ as an abstract group consists of the defining relations of $A$ and $B$,
together with the relations which glue together the copies in $A$ and $B$
of the elements of $C$. This expresses the defining relations of $\text{Aut }K[x,y]$
in terms of the defining relations of the groups of the affine and the triangular
automorphisms. The picture is much more complicated in the case of tame automorphisms of
$K[x,y,z]$. Analyzing the algorithm of Theorem \ref{theorem of Shestakov and Umirbaev},
very recently Umirbaev \cite{U2} obtained a set of defining relations of $\text{TAut }K[x,y,z]$.

We write $x_1,x_2,x_3$ for the variables $x,y,z$ and define the automorphisms
\begin{equation}\label{generators of the tame group}
\sigma(i,\alpha,f)=(x_1,\ldots,\alpha x_i+f,\ldots,x_3)
\end{equation}
of the polynomial algebra $K[x_1,x_2,x_3]$,
where $i=1,2,3$, $\alpha\in K^{\ast}$, and $f\in K[x_1,x_2,x_3]$ does not depend on
the variable $x_i$. Clearly, the automorphisms
$\sigma(i,\alpha,f)$ generate $\text{TAut }K[x_1,x_2,x_3]$.
Let us denote
\[
\tau_{(ks)}=\sigma(s,-1,x_k)\sigma(k,1,-x_s)\sigma(s,1,x_k),\quad k,s=1,2,3,\quad k\not=s.
\]
The only nontrivial action of $\tau_{(ks)}$ is to change the variables $x_k$ and $x_s$.
The defining relations of $\text{TAut }K[x_1,x_2,x_3]$ are given by
the following important theorem.

\begin{theorem}\label{defining relations for tame automorphisms}
{\rm (Umirbaev \cite{U2})}
Let $K$ be a field of characteristic $0$.
The group of tame automorphisms of $K[x_1,x_2,x_3]$ has the following set of defining
relations with respect to the set of generators (\ref{generators of the tame group}):
\begin{equation}\label{first defining relation}
\sigma(i,\alpha,f)\sigma(i,\beta,g)=\sigma(i,\alpha\beta,\beta f+g);
\end{equation}
If $i\not=j$ and $f$ does not depend on $x_i,x_j$, then
\begin{equation}\label{second defining relation}
\sigma(i,\alpha,f)^{-1}\sigma(j,\beta,g)\sigma(i,\alpha,f)=
\sigma(j,\beta,\sigma(i,\alpha,f)^{-1}(g));
\end{equation}
If $x_j=\tau(ks)(x_i)$, then
\begin{equation}\label{third defining relation}
\tau_{(ks)}\sigma(i,\alpha,f)\tau_{(ks)}=\sigma(j,\alpha,\tau_{(ks)}(f)).
\end{equation}
\end{theorem}

\section{A survey on automorphisms of free algebras}

Let $K\langle x_1,\ldots,x_n\rangle$ be the free associative algebra freely generated
by $x_1,\ldots,x_n$ over an arbitrary field $K$.
We may think of it as the algebra of polynomials in $n$ noncommuting variables.
As in the commutative case,
the group $\text{TAut }K\langle x_1,\ldots,x_n\rangle$ of the tame automorphisms
is generated by the affine and the triangular automorphisms.
We start the section with the description of the automorphisms and coordinates
of $K\langle x,y\rangle$. Recall that in noncommutative algebra
coordinates are also called primitive elements.

\begin{theorem}\label{theorem of Czerniakiewicz and Makar-Limanov}
{\rm (Czerniakiewicz \cite{Cz} and Makar-Limanov \cite{ML1, ML2})}
The automorphisms of $K\langle x,y\rangle$ are tame. The groups
$\text{\rm Aut }K\langle x,y\rangle$ and $\text{\rm Aut }K[x,y]$ are isomorphic.
\end{theorem}

The idea of the proof is the following.
The natural homomorphism
\begin{equation}\label{abelianization}
\pi: K\langle x,y\rangle\to K[x,y]
\end{equation}
induces a group homomorphism
\begin{equation}\label{isomorphism of automorphism groups}
\pi_1:\text{Aut }K\langle x,y\rangle\to \text{Aut }K[x,y]
\end{equation}
Every automorphism of $K[x,y]$ is tame and hence can be lifted to a tame
automorphism of $K\langle x,y\rangle$. Hence the mapping $\pi_1$ is onto.
It is sufficiently to show that $\text{Ker }\pi_1=1$.
The kernel of $\pi$ is the commutator ideal of $K\langle x,y\rangle$,
generated as an ideal by the commutator $[x,y]=xy-yx$.
Hence one has to show that there are no automorphisms of
$K\langle x,y\rangle$ of the form
$\varphi=(x+u,y+v)$, where $u,v$ are in the commutator ideal
and at least one of them is different from 0. The proof is based
on combinatorial analysis of the words in the free algebra
and uses essentially the weak Euclidean algorithm. For details we refer to
Chapter 6 of the book by Cohn \cite{C2}.

The following theorem provides a simple test whether an endomorphism
of $K\langle x,y\rangle$ is an automorphism.

\begin{theorem}\label{commutator test}
{\rm (Dicks \cite{Di})}
The endomorphism $\varphi$ of $K\langle x,y\rangle$ is an automorphism
if and only if there exists a constant $\alpha\in K^{\ast}$ such that
\[
[\varphi(x),\varphi(y)]=\alpha[x,y].
\]
\end{theorem}

Unfortunately the result does not give any decomposition
of the automorphism as a product of
affine and triangular automorphisms.

Algorithms \ref{algorithm for recognizing of automorphisms} and
\ref{algorithm for recognizing of coordinates} can be easily modified to
recognize the automorphisms and the coordinates of $K\langle x,y\rangle$.
The only difference is that one has to work in noncommutative
setup. For example, in the first algorithm we have to determine whether
for two homogeneous polynomials $u,v\in K\langle x,y\rangle$ one of them
is a power of the other, e.g. $u=\beta v^d$ for some $\beta\in K^{\ast}$.
Since $d$ is determined by the degrees of $u$ and $v$, it is sufficient to write
both $u$ and $v^d$ as elements in $K\langle x,y\rangle$ and to see whether
the corresponding monomials are proportional. In the second algorithm
one has to check whether some $(d,1)$-homogeneous polynomial is of the form
$\eta(x+\xi y^d)$ and this also can be done easily in $K\langle x,y\rangle$.

An easier way to handle problems for the automorphisms and coordinates of
$K\langle x,y\rangle$ is to reduce the considerations to similar problems in
$K[x,y]$ and then to use the isomorphism
(\ref{isomorphism of automorphism groups}) of the automorphism groups of
$K[x,y]$ and $K\langle x,y\rangle$. This approach was applied by Shpilrain and Yu \cite{SY2}.
In particular, they found the first algorithm which recognizes the coordinates
of $K\langle x,y\rangle$.

\begin{algorithm}\label{algorithm for automorphisms of free algebras of rank 2}
Let $\varphi=(f,g)$ be an endomorphism of the free algebra
$K\langle x,y\rangle$.

{\it Step 1}. Consider the endomorphism
$\psi=(\pi(f),\pi(g))$ of $K[x,y]$, where $\pi(f)$ and $\pi(g)$ are
the abelianization of $f$ and $g$ from (\ref{abelianization}).
Apply Algorithm \ref{algorithm for recognizing of automorphisms} to $\psi$.
If $\psi$ is not an automorphism of $K[x,y]$, then
$\varphi$ is not an automorphism of $K\langle x,y\rangle$.

{\it Step 2}. If $\psi$ is an automorphism, decompose it as a product
$\psi=\psi_1\cdots\psi_n$ of affine and triangular automorphisms $\psi_i$.
Consider the automorphism
\[
\pi_1^{-1}(\psi)=\pi_1^{-1}(\psi_1)\cdots\pi_1^{-1}(\psi_n).
\]
Then $\varphi$ is an automorphism if and only if
$\varphi=\pi_1^{-1}(\psi)$.
\end{algorithm}

\begin{algorithm}\label{algorithm for coordinates of free algebras of rank 2}
(Shpilrain and Yu \cite{SY2})
Let $f\in K\langle x,y\rangle$.

{\it Step 1}. Consider the abelianization $\pi(f)$ of $f$.
Apply Algorithm \ref{algorithm for recognizing of coordinates},
or \ref{algorithm for coordinates in two variables}
if $\text{char }K=0$,
to decide whether $\pi(f)$ is a coordinate of $K[x,y]$.
If it is not, then $f$ is not a coordinate of $K\langle x,y\rangle$.

{\it Step 2}. If $\pi(f)$ is a coordinate of $K[x,y]$, find an
automorphism $\psi$ of $K[x,y]$ which sends $x$ to $\pi(f)$.
Then $f$ is a coordinate if and only if $\varphi=\pi_1^{-1}(\psi)$
sends $x$ to $f$.
\end{algorithm}

For the motivation of the algorithm,
let $\varphi=\pi_1^{-1}(\psi)\in\text{Aut }K\langle x,y\rangle$, where
$\psi\in \text{Aut }K[x,y]$ sends $x$ to $\pi(f)$.
All automorphisms $\psi'\in \text{Aut }K[x,y]$ with this property are of the form
$\psi'=\psi\cdot(x,\beta y+h(x))$, $\beta\in K^{\ast}$, $h(x)\in K[x]$.
Hence
\[
\varphi'=\pi_1^{-1}(\psi')=\pi_1^{-1}(\psi)\cdot(x,\beta y+h(x))
\varphi\cdot(x,\beta y+h(x))
\]
and $\varphi(x)=\varphi'(x)$. See \cite{SY2} also for other properties of the coordinates
of $K\langle x,y\rangle$.

As in the commutative case, see \cite{Wi, MLSY}, one considers the problem whether
two polynomials $p,q\in K\langle x,y\rangle$ are automorphically equivalent.
The corresponding algorithm is given by Drensky and Yu \cite{DY6}.

By analogy with the commutative case one can ask immediately
the following problem. It is stated, for example, by Cohn
as Problem 1.74 in \cite{DN}.

\begin{problem}
Is every automorphism of $K\langle x_1,\ldots,x_n\rangle$, $n>2$, tame?
\end{problem}

If we are able to lift any wild automorphism of $K[x,y,z]$ to an automorphism
of $K\langle x,y,z\rangle$, then we shall obtain automatically
a wild automorphism. Up till now, no such wild automorphism has been
constructed.

\begin{problem}\label{can one lift Nagata automorphism}
Can the wild automorphisms of $K[x,y,z]$ be lifted to automorphisms of
$K\langle x,y,z\rangle$? Can the Nagata automorphism be lifted?
Can the wild $z$-automorphisms of $K[z][x,y]$ be lifted to $z$-automorphisms
of $K\langle x,y,z\rangle$ (or at least to any automorphisms
of $K\langle x,y,z\rangle$, not necessarily fixing $z$)?
\end{problem}

A stronger version of this problem is the following.

\begin{problem}
Can every wild coordinate of $K[x,y,z]$ be lifted to a coordinate of
$K\langle x,y,z\rangle$? Can the two nontrivial coordinates of the Nagata automorphism be lifted?
Can the two nontrivial coordinates of every wild $z$-automorphism of $K[z][x,y]$ be lifted to
coordinates of a $z$-automorphism of $K\langle x,y,z\rangle$?
\end{problem}

The most prospective candidate for a wild automorphism
of the free algebra with three generators is the example of Anick
$(x+(y(xy-yz),y,z+(zy-yz)y)$,
see the book by Cohn \cite{C2}, p.~343, or,
$(x+z(xz-zy),y+(xz-zy)z,z)$,
changing the places of $y$ and $z$.

\begin{conjecture}\label{Anick conjecture}
The Anick automorphism
\[
\nu_1=(x+z(xz-zy),y+(xz-zy)z,z)
\]
of $K\langle x,y,z\rangle$ is wild.
\end{conjecture}

The Anick automorphism fixes $z$. Hence we may consider
the variable $z$ as a noncommutative constant and treat the Anick
automorphism as a $z$-automorphism
of $K\langle x,y,z\rangle$. It has the property that the first two coordinates
are linear in $x$ and $y$. We call such automorphisms linear $z$-automorphisms. From
some point of view, the linear $z$-automorphisms
are the simplest $z$-automorphisms. It is natural to expect the appearance of
some group of matrices there.

There is a general concept of partial derivatives and the Jacobian matrix
in the free algebra, developed by Yagzhev \cite{Y},
Dicks and Lewin \cite{DiL}, and Schofield \cite{Sc}, including the solution of
the Jacobian conjecture for free associative algebras for $n=2$ in \cite{DiL}
and arbitrary $n$ in \cite{Sc}.

For the study of linear $z$-automorphisms,
Drensky and Yu \cite{DY3} introduced a simplified version of the
partial derivatives of Dicks and Lewin \cite{DiL}.
We shall restrict our considerations to the free algebra with three generators only.
Let $f\in K\langle x,y,z\rangle$ be linear in $x,y$. Then it has the form
\begin{equation}\label{linear z-polynomial}
f=\sum\alpha_{ij}z^ixz^j+\sum\beta_{ij}z^iyz^j,\quad \alpha_{ij},\beta_{ij}\in K.
\end{equation}
The $z$-derivatives
$f_x$ and $f_y$ are defined by
\begin{equation}\label{partial z-derivatives}
f_x=\sum\alpha_{ij}z_1^iz_2^j,\quad
f_y=\sum\beta_{ij}z_1^iz_2^j.
\end{equation}
Here $f_x$ and $f_y$ are in $K[z_1,z_2]$ and are polynomials in two commuting variables.
The $z$-Jacobian matrix of the linear $z$-endomorphism $\varphi=(f,g)$
of $K\langle x,y,z\rangle$ is defined as
\begin{equation}\label{z-Jacobian}
J_z(\varphi)=\left(
\begin{matrix}
f_x&g_x\\
f_y&g_y\\
\end{matrix}
\right).
\end{equation}

It is easy to see that the linear $z$-endomorphism $\varphi$ of $K\langle x,y,z\rangle$
is an automorphism if and only if the matrix $J_z(\varphi)$ is invertible as a matrix
with entries in $K[z_1,z_2]$. Also, the group of linear $z$-automorphisms is
isomorphic to $GL_2(K[z_1,z_2])$. The following theorem was established in \cite{DY3}.

\begin{theorem}\label{linear z-wild automorphisms}
{\rm (i)} A linear $z$-automorphism is $z$-tame if and only if
its $z$-Jacobian matrix belongs to $GE_2(K[z_1,z_2])$.

{\rm (ii)} Every linear $z$-automorphism of $K\langle x,y,z\rangle$
induces a tame automorphism of $K[x,y,z]$.
\end{theorem}

Applying the result of
Tolhuizen, Hollmann and Kalker \cite{THK}
and Park \cite{P} we can decide whether
a linear $z$-automorphism is $z$-tame applying
the Euclidean division algorithm to the $z$-Jacobian matrix.
We apply this to the Anick automorphism.

\begin{corollary}\label{Anick automorphism is z-wild}
The Anick automorphism is $z$-wild.
\end{corollary}

Really, the $z$-Jacobian matrix of the Anick automorphism is
\[
J_z(\nu_1)=\left(
\begin{matrix}
1+z_1z_2&z_2^2\\
-z_1^2&1-z_1z_2\\
\end{matrix}
\right).
\]
The matrix $J_z(\nu_1)$ coincides with the well known example of Cohn
\cite{C1}
of a matrix in $GL_2(K[z_1,z_2])$ which does not belong to $GE_2(K[z_1,z_2])$.
(Direct arguments: Since the homogeneous components of second (maximal) degree of
the first column of $J_z(\nu_1)$ do not divide each other, we cannot
apply the Euclidean algorithm.) Hence the Anick automorphism is $z$-wild.
The same arguments work for the matrix
\[
J_z(\varphi)=\left(
\begin{matrix}
1+z_1z_2h(z_1,z_2)&z_2^2h(z_1,z_2)\\
-z_1^2h(z_1,z_2)&1-z_1z_2h(z_1,z_2)\\
\end{matrix}
\right),
\]
where $h(z_1,z_2)$ is any nonzero polynomial in $K[z_1,z_2]$. Then
the corresponding automorphism $\varphi$ is $z$-wild.

By the celebrated theorem of Suslin \cite{Su}, every matrix in
$GL_n(K[z_1,\ldots,z_p])$, $n\geq 3$,
can be presented as a product of a diagonal matrix
and elementary matrices. Hence, for $n\geq 3$ all
linear $z$-automorphisms of $K\langle x_1,\ldots,x_n,z\rangle$ are tame.
In particular, the linear $z$-automorphisms of $K\langle x,y,z\rangle$ are stably tame.
For example, from a lemma of Mennicke, used in the constructive proof of the Suslin theorem
given by Park and Woodburn \cite{PW}, follows that, extended by $\nu_1(t)=t$,
the Anick automorphism of $K\langle x,y,t,z\rangle$ can be decomposed as a
product of ``elementary'' $z$-automorphisms
\[
\nu_1=
\varepsilon_{31}(z_1)\varepsilon_{32}(z_2)\varepsilon_{13}(-z_2)
\varepsilon_{23}(z_1)\varepsilon_{31}(-z_1)\varepsilon_{32}(-z_2)
\varepsilon_{13}(z_2)\varepsilon_{23}(-z_1),
\]
where the the $z$-automorphism $\varepsilon_{ij}(\alpha z_1^az_2^b)$
of $K\langle x,y,t,z\rangle$ is defined by
\[
\varepsilon_{ij}(\alpha z_1^az_2^b):x_j\to x_j+\alpha z^ax_iz^b,
\quad  \varepsilon_{ij}(\alpha z_1^az_2^b): x_k\to x_k, \quad k\not=j,
\]
and $x_1=x$, $x_2=y$, $x_3=t$.

The proof of Theorem \ref{linear z-wild automorphisms}
in \cite{DY3}
gives a criterion for $z$-wildness of a class of nonlinear $z$-automorphisms.
See the discussions below on our results in
\cite{DY4, DY5} for further development of the idea.

\begin{corollary}\label{z-tame implies linear z-tame}
Let $\varphi=(f,g)$ be a $z$-automorphism of $K\langle x,y,z\rangle$. Let
$\varphi_1=(f_1,g_1)$, where $f_1,g_1$ are the linear in $x$ and $y$
components of $f,g$. If the matrix $J_z(\varphi_1)$ does not belong to
$GE_2(K[z_1,z_2])$, then the automorphism $\varphi$ is $z$-wild.
\end{corollary}

The study of $z$-tame automorphisms of $K\langle x,y,z\rangle$ has continued in the
very recent paper by Drensky and Yu \cite{DY7}. The following result is an analogue of
Theorem \ref{theorem of Wright} (i) established by Wright \cite{Wr}.
It shows that the structure of the group
of $z$-tame automorphisms of $K\langle x,y,z\rangle$ is similar to the structure of
the group of $z$-tame automorphisms of $K[x,y,z]$.

\begin{theorem}\label{structure of z-TAut} {\rm (Drensky and Yu \cite{DY7})}
Over an arbitrary field $K$,
the group $\text{\rm TAut}_zK\langle x,y,z\rangle$ of $z$-tame
automorphisms of $K\langle x,y,z\rangle$
is isomorphic to the free product $A\ast_CB$ of the group $A$
of the $z$-affine automorphisms and the group $B$ of $z$-triangular automorphisms
amalgamating their intersection $C=A\cap B$.
\end{theorem}

Since the group of the $z$-automorphisms which
are linear in $x,y$ is isomorphic to the group $GL_2(K[z_1,z_2])$,
we immediately obtain:

\begin{corollary}
The group $\text{\rm TAut}_zK\langle x,y,z\rangle$
is isomorphic to the free product with amalgamation
$GE_2(K[z_1,z_2])\ast_{C_1}B$, where
$GE_2(K[z_1,z_2])$ is identified as above with the group of $z$-tame automorphisms
which are linear in $x$ and $y$, $B$ is the group of $z$-triangular automorphisms and
$C_1=GE_2(K[z_1,z_2])\cap B$.
\end{corollary}

Now we can use Theorem \ref{structure of z-TAut} to present algorithms
which recognize $z$-tame automorphisms and
coordinates of $K\langle x,y,z\rangle$.
We start with an algorithm
which determines whether a $z$-endomorphism of $K\langle x,y,z\rangle$
is a $z$-tame automorphism. The main idea is similar to that of
Algorithm \ref{algorithm for recognizing of automorphisms}
which decides whether an endomorphism of
$K[x,y]$ is an automorphism, but the realization
is more sophisticated. Instead of the usual degree used there
we define a bidegree of $K\langle x,y,z\rangle$ assuming that the monomial $w$ is
of bidegree $\text{bideg }w=(d,e)$ if
$\text{deg}_xw+\text{deg}_yw=d$ and $\text{deg}_zw=e$.
We order the bidegrees $(d,e)$ lexicographically, i.e.,
$(d_1,e_1)>(d_2,e_2)$ means that either $d_1>d_2$ or
$d_1=d_2$ and $e_1>e_2$.
In order to simplify further the considerations, we
use the trick introduced by Formanek \cite{Fo} in his
construction of central polynomials of matrices.

Let $H_n$ be the subspace of $K\langle x,y,z\rangle$ consisting of all
polynomials which are homogeneous of degree $n$ with respect to $x$ and $y$.
We define an action of $K[t_0,t_1,\ldots,t_n]$ on $H_n$ in the following way. If
\[
w=z^{a_0}u_1z^{a_1}u_2\cdots z^{a_{n-1}}u_nz^{a_n},
\]
where $u_i=x$ or $u_i=y$, $i=1,\ldots,n$, then
\[
t_0^{b_0}t_1^{b_1}\cdots t_n^{b_n}\ast w=
z^{a_0+b_0}u_1z^{a_1+b_1}u_2\cdots z^{a_{n-1}+b_{n-1}}u_nz^{a_n+b_n},
\]
and extend this action by linearity. Clearly, $H_n$ is a free
$K[t_0,t_1,\ldots,t_n]$-module with basis consisting of the $2^n$
monomials $u_1\cdots u_n$, where $u_i=x$ or $u_i=y$.
The proofs of the following equations are obtained by easy direct computation.
Let $\beta\in K^{\ast}$,
\begin{equation}\label{this is v}
v(x,y,z)=\sum\theta_i(t_0,t_1,\ldots,t_k)\ast u_{i_1}\cdots u_{i_k}\in H_k,
\end{equation}
\begin{equation}\label{this is q}
q(y,z)=\omega(t_0,t_1,\ldots,t_d)\ast y^d\in H_d,
\end{equation}
where $\theta_i\in K[t_0,t_1,\ldots,t_k]$, $\omega\in K[t_0,t_1,\ldots,t_d]$,
$u_{i_j}=x$ or $u_{i_j}=y$. Then
\begin{equation}\label{this is u}
\begin{array}{c}
u(x,y,z)=q(v(x,y,z)/\beta,z)
=\omega(t_0,t_d,t_{2d},\ldots,t_{kd})/\beta^d\\
\\
\left(\sum\theta_i(t_0,t_1,\ldots,t_k)\ast u_{i_1}\cdots u_{i_k}\right)\\
\\
\left(\sum\theta_i(t_k,t_{k+1},\ldots,t_{2k})\ast u_{i_1}\cdots u_{i_k}\right)\cdots\\
\\
\left(\sum\theta_i(t_{k(d-1)},t_{k(d-1)+1},\ldots,t_{kd})\ast u_{i_1}\cdots u_{i_k}\right).\\
\end{array}
\end{equation}

\begin{algorithm}\label{algorithm for z-tame automorphisms} {\rm (Drensky and Yu \cite{DY7})}
Let $\varphi=(f,g)$ be a $z$-endomorphism of $K\langle x,y,z\rangle$.
We make use of the bidegree defined above.

{\it Step 0}. If some of the polynomials $f,g$ depends on $z$ only,
then $\varphi$ is not an automorphism.

{\it Step 1}. Let $u,v$ be the homogeneous components of highest bidegree of $f,g$,
respectively. If both $u,v$ are of bidegree $(1,0)$, i.e., linear, then we check
whether they are linearly
independent. If yes, then $\varphi$ is a product of a linear automorphism
(from $GL_2(K)$) and a translation $(x+p(z),y+r(z))$. If
$u,v$ are linearly dependent, then $\varphi$ is not an automorphism.

{\it Step 2}. Let $\text{bideg }u>(1,0)$ and
$\text{bideg }u\geq \text{bideg }v$. Hence $u\in H_l$, $v\in H_k$ for some $k$ and $l$.
We have to check whether $l=kd$ for a positive integer $d$
and to decide whether $u=q(v/\beta,z)$ for some
$\beta\in K^{\ast}$ and some $q(y,z)\in H_d$. We know
$u$ in (\ref{this is u}) and $v$ in (\ref{this is v}) up to the multiplicative constant $\beta$.
Hence, up to $\beta$, we know the polynomials
$\theta_i(t_0,t_1,\ldots,t_n)$ in the presentation of $v$.
We compare some of the nonzero polynomial coefficients of
$u=\sum\lambda_j(t_0,\ldots,t_{kd})u_{j_1}\cdots u_{i_{kd}}$
with the corresponding coefficient of $q(v/\beta,z)$.
We can find explicitly, up to
the value of $\beta^d$, the polynomial
$\omega(t_0,t_1,\ldots,\omega_d)$ in (\ref{this is q}) using the usual division of polynomials.
If $l=kd$ and $u=q(v/\beta,z)$, then we replace
$\varphi=(f,g)$ with $\varphi_1=(f-q(g/\beta,z),g)$.
Then we apply Step 0 to $\varphi_1$.
If $u$ cannot be presented in the desired form,
then $\varphi$ is not an automorphism.

{\it Step 3}. If $\text{bideg }v>(1,0)$ and
$\text{bideg }u<\text{bideg }v$, we have similar considerations, as
in Step 2, replacing
$\varphi=(f,g)$ with $\varphi_1=(f,g-q(f/\alpha,z))$ for a suitable
$q(y,z)$.
Then we apply Step 0 to $\varphi_1$.
If $v$ cannot be presented in this form,
then $\varphi$ is not an automorphism.
\end{algorithm}

The following example from \cite{DY4, DY5, DY7} is based on the polynomial $xz-zy$ which
appears in the Anick automorphism.

\begin{example}\label{nonmetabelian example generalizing Anick}
Let $h(t,z)\in K\langle t,z\rangle$.
Then
\[
\sigma_h=(x+zh(xz-zy,z),y+h(xz-zy,z)z,z)
\]
is an automorphism of $K\langle x,y,z\rangle$  fixing $xz-zy$.
Its inverse is $\sigma_{-h}$.
\end{example}

\begin{corollary}\label{nonmetabelian generalization of Anick}
{(\rm Drensky and Yu \cite{DY7})}
Let $h(t,z)\in K\langle t,z\rangle$ and let $\text{\rm deg}_uh(u,z)>0$.
Then
\[
\sigma_h=(x+zh(xz-zy,z),y+h(xz-zy,z)z,z)
\]
is a $z$-wild automorphism of $K\langle x,y,z\rangle$.
\end{corollary}

For the proof,
we apply Algorithm \ref{algorithm for z-tame automorphisms}. Let $w$ be the
homogeneous component of highest bidegree of $h(xz-zy,z)$. Clearly,
$w$ has the form $w=\overline{h}(xz-zy,z)=q(xz-zy,z)$ for some bihomogeneous
polynomial $q(t,z)\in K\langle t,z\rangle$. The leading components
of the coordinates of $\sigma_h$ are $zq(xz-zy,z)$ and $q(xz-zy,z)z$, and are of
the same bidegree. If $\sigma_h$ is a $z$-tame automorphism, then we can
reduce the bidegree using a linear transformation, which is impossible because
$zq(xz-zy,z)$ and $q(xz-zy,z)z$ are linearly independent.

Clearly, if $\varphi=(f,g)$ is a $z$-automorphism
of $K\langle x,y,z\rangle$, with the
same first coordinate, then we can find all other $z$-coordinate mates of $f$.
This argument and Corollary \ref{nonmetabelian generalization of Anick} give immediately:

\begin{corollary}\label{coordinates of nonmetabelian generalization of Anick}
{(\rm Drensky and Yu \cite{DY7})}
Let $h(t,z)\in K\langle t,z\rangle$ and let $\text{\rm deg}_uh(u,z)>0$.
Then $f(x,y,z)=x+zh(xz-zy,z)$
is a $z$-wild coordinate of $K\langle x,y,z\rangle$.
\end{corollary}

Now we are able to modify
Algorithm \ref{algorithm for z-tame automorphisms}
to decide whether a polynomial $f(x,y,z)$ is a $z$-tame coordinate
of $K\langle x,y,z\rangle$. We refer to \cite{DY7} for details.

\begin{theorem} {(\rm Drensky and Yu \cite{DY7})}
There is an algorithm which decides whether a polynomial
$f(x,y,z)\in K\langle x,y,z\rangle$ is a $z$-tame coordinate.
\end{theorem}

As we already mentioned in the comments after Theorem \ref{theorem of Czerniakiewicz and Makar-Limanov},
the main step of the proof of
Czerniakiewicz \cite{Cz} and Makar-Limanov \cite{ML1, ML2}
for the tameness of the automorphisms of $K\langle x,y\rangle$ is the following.
If $\varphi=(x+u,y+v)$ is an endomorphism of $K\langle x,y\rangle$,
where $u,v$ are in the commutator ideal of $K\langle x,y\rangle$
and at least one of them is different from 0,
then $\varphi$ is not an automorphism of $K\langle x,y\rangle$.
The condition that $u(x,y)$ and $v(x,y)$ belong to the commutator ideal of $K\langle x,y\rangle$
immediately implies that
all monomials of $u$ and $v$ depend on both $x$ and $y$, as required in the
following theorem.

\begin{theorem}
{(\rm Drensky and Yu \cite{DY7})}
The $z$-endomorphisms of the form
\[
\varphi=(x+u(x,y,z),y+v(x,y,z)),
\]
where $(u,v)\not=(0,0)$ and all monomials of $u$ and $v$ depend on both $x$ and $y$, are not automorphisms
of $K\langle x,y,z\rangle$.
\end{theorem}

Till the end of the section we fix the field $K$ of characteristic 0
and continue our survey with the work of Umirbaev \cite{U2} devoted on the wildness
of the Anick automorphism. It is too difficult to controll the behavior of
the tame automorphisms of $K\langle x,y,z\rangle$. Instead, it is more convenient
to consider some factor algebra of $K\langle x,y,z\rangle$ modulo an ideal
of $K\langle x,y,z\rangle$ which is invariant under the action of the automorphism
group of $K\langle x,y,z\rangle$.
This ideal should be big enough to allow to work
in the factor algebra, and not too big,
in order to preserve the wildness of some automorphisms.
The ideals of the free algebra which are invariant
under all endomorphisms are called T-ideals.
Such an ideal $U$ of $K\langle x_1,\ldots,x_n\rangle$
can be characterized with the property that there exists an associative algebra
$R$ such that $U$ coincides with the ideal of all polynomial polynomial identities
$u(x_1,\ldots,x_n)=0$ of $R$. This means that $u(r_1,\ldots,r_n)=0$
for all $r_1,\ldots,r_n\in R$. The factor algebra
$K\langle x_1,\ldots,x_n\rangle/U$ is called the (relatively) free algebra of rank $n$
in the variety of algebras generated by the algebra $R$, or defined by the polynomial
identities from $U$.

The obtained results have shown that it is very convenient to work in the free metabelian
associative algebra which we denote by $C_n$.
The variety of metabelian associative algebras is determined
by the polynomial identity $[t_1,t_2][t_3,t_4]=0$.
One can define partial derivatives and the related Jacobian matrix
induced by the corresponding objects in the free associative algebra.
Although $C_n$ is a homomorphic image of $K\langle x_1,\ldots,x_n\rangle$, we denote its
generators with the same symbols $x_1,\ldots,x_n$.
For the partial derivatives we need two sets
of commuting variables $\{u_1,\ldots,u_n\}$ and $\{v_1,\ldots,v_n\}$.
We define formal partial derivatives
$\partial_M/\partial_Mx_i$ by
\[
\frac{\partial_Mx_i}{\partial_Mx_i}=1,
\quad \frac{\partial_Mx_j}{\partial_Mx_i}=0,\,j\not=i,
\]
and, for  a monomial
$w=x_{i_1}\cdots x_{i_m}\in C_n$,
\[
\frac{\partial_Mw}{\partial_Mx_i}=
\sum_{k=1}^mu_{i_1}\cdots u_{i_{k-1}}
v_{i_{k+1}}\ldots v_{i_m}
\frac{\partial_Mx_{i_k}}{\partial_Mx_i}.
\]
A polynomial $f\in C_n$ belongs to the commutator ideal of $C_n$,
i.e., to the kernel of the natural homomorphism
$C_n\to K[x_1,\ldots,x_n]$, if and only if
\[
\sum_{i=1}^n(u_i-v_i)\frac{\partial_Mf}{\partial_Mx_i}=0.
\]
The Jacobian matrix of an endomorphism $\varphi$ of $C_n$ is
\begin{equation}\label{metabelian Jacobian}
J_M(\varphi)=\left(\frac{\partial_M\varphi(x_j)}{\partial_Mx_i}\right),
\end{equation}
which is a matrix with entries from the polynomial algebra
$K[u_1,\ldots,u_n,v_1,\ldots,v_n]$.
Umirbaev \cite{U1} proved that the
Jacobian matrix $J_M(\varphi)$ is invertible (as a matrix
with polynomial entries) if and only if
$\varphi$ is an automorphism of $C_n$. Also,
every automorphism of $K[x_1,\ldots,x_n]$ can be lifted to an automorphism
of $C_n$. Clearly, the invertibility of $J_M(\varphi)$ is equivalent
to $0\not= \text{\rm det}(J_M(\varphi))\in K$.

We shall work with free algebras of rank 3 only and
shall use the sets of variables
\[
X=\{x,y,z\},\quad U=\{x_1,y_1,z_1\},\quad V=\{x_2,y_2,z_2\}.
\]
instead of $X=\{x_1,x_2,x_3\}$, $U=\{u_1,u_2,u_3\}$, $V=\{v_1,v_2,v_3\}$,
respectively, and shall denote the algebra $C_3$ by $C$.
For example, if
$f\in C$ is linear in $x$ and $y$, it can be written in the form
(\ref{linear z-polynomial}) and its partial derivatives
$\partial_Mf/\partial_Mx$ and $\partial_Mf/\partial_My$ coincide, respectively,
with $f_x$ and $f_y$ in (\ref{partial z-derivatives}). Hence, if
$\varphi=(f,g,z)$ is a $z$-endomorphism of $C$, and $f,g$ are linear
in $x,y$, then the Jacobian matrix $J_M(\varphi)$ is of the form
\begin{equation}\label{metabelian z-Jacobian}
J_M(\varphi)=\left(
\begin{matrix}
f_x&g_x&0\\
f_y&g_y&0\\
\partial_Mf/\partial_Mz&\partial_Mg/\partial_Mz&1\\
\end{matrix}
\right),
\end{equation}
compare it with $J_z(\varphi)$ from (\ref{z-Jacobian}).

There is a natural homomorphism
\[
\text{TAut }K\langle x,y,z\rangle\to \text{TAut}(C)\to \text{TAut }K[x,y,z].
\]
As a consequence of Theorem \ref{defining relations for tame automorphisms},
Umirbaev \cite{U2} obtained the following.

\begin{proposition}\label{the kernel of tame automorphisms}
As a normal subgroup of $\text{\rm TAut}(C)$, the kernel of
the homomorphism
\[
\text{\rm TAut}(C)\to \text{\rm TAut }K[x,y,z]
\]
is generated by the automorphisms
\[
\psi_w=(x+w(y,z),y,z),\quad w(y,z)
=\sum_{p,q,r,s\geq 0}\alpha_{pqrs}y^pz^q[y,z]y^rz^s,\,\alpha_{pqrs}\in
K.
\]
\end{proposition}

Hence the kernel of the homomorphism
$\text{\rm TAut }(C)\to \text{\rm TAut }K[x,y,z]$ is generated by
automorphisms of the form
\begin{equation}\label{typical generator of the kernel}
\vartheta=\sigma(i_k,\alpha_k,g_k)^{-1}\cdots\sigma(i_1,\alpha_1,g_1)^{-1}\psi_w
\sigma(i_1,\alpha_1,g_1)\cdots\sigma(i_k,\alpha_k,g_k)
\end{equation}
for some automorphisms $\sigma(i_j,\alpha_j,g_j)$
from (\ref{generators of the tame group}).
The Jacobian matrix of $\psi_f$ is of the form
\[
J_M(\psi_w)=\left(
\begin{matrix}
1&0&0\\
\partial_Mw/\partial_My&1&0\\
\partial_Mw/\partial_Mz&0&1\\
\end{matrix}
\right),
\]
which $\sigma(i_j,\alpha_j,g_j)$ also has an elementary Jacobian matrix.
Unfortunately, this fact cannot be applied directly to recognize the tame automorphisms
due to the validity of the Suslin theorem \cite{Su} for $n\geq 3$.
On the other hand, in the case of the Anick automorphism and other linear $z$-automorphisms
$\varphi$ the relation between $J_z(\varphi)$ from (\ref{z-Jacobian})
and $J_M(\varphi)$ from (\ref{metabelian Jacobian}) is obvious. Hence
it is natural to try to apply the idea of Theorem \ref{linear z-wild automorphisms}
and to reduce the considerations to the $2\times 2$ case.
This was done by Umirbaev \cite{U2} in the following way.
For an automorphism $\varphi=(f,g,h)$ of $C$, he introduced
the $2\times 2$ ``Jacobian'' matrix
\begin{equation}\label{Jacobian of Umirbaev}
J_2(\varphi)=\left(
\begin{matrix}
\partial_Mf/\partial_Mx&\partial_Mg/\partial_Mx\\
\partial_Mf/\partial_My&\partial_Mg/\partial_My\\
\end{matrix}
\right).
\end{equation}
Then he defined a homomorphism $\eta:K[z_1,y_1,z_1,x_2,y_2,z_2]\to K[z_1,z_2]$
by
\[
\eta(x_1)=\eta(y_1)=\eta(x_2)=\eta(y_2)=0,\quad \eta(z_1)=z_1,\quad \eta(z_2)=z_2.
\]
Using the concrete form of the generators (\ref{generators of the tame group})
of $\text{TAut }K[x,y,z]$ and the defining relations
(\ref{first defining relation}), (\ref{second defining relation}),
and (\ref{third defining relation}) between them, Umirbaev established that
for the automorphism $\vartheta$ from (\ref{typical generator of the kernel})
the matrix $\eta\left(J_2(\vartheta)\right)$ belongs to $GE_2(K[z_1,z_2])$.
This implies the following theorem.

\begin{theorem}\label{criterion for tameness in metabelian case}
{\rm (Umirbaev \cite{U2})}
Let the field $K$ be of characteristic $0$
and let $\varphi$ be a tame automorphism of the free metabelian algebra $C$
which induces the identity automorphism of $K[x,y,z]$. Then the
matrix $\eta\left(J_2(\varphi)\right)$
belongs to the group $GE_2(K[z_1,z_2])$.
\end{theorem}

Now, there is only one obvious step left to the proof of the wildness of the Anick automorphism.
Consider a tame linear $z$-automorphism $\vartheta$ of $K\langle x,y,z\rangle$
which induces the same automorphism of $K[x,y,z]$ as the Anick automorphism $\nu_1$.
One can apply Theorem \ref{criterion for tameness in metabelian case} to
$\vartheta^{-1}\nu_1$. But the matrix
\[
\eta\left(J_2(\vartheta^{-1}\nu_1)\right)=
\eta\left(J_2(\vartheta^{-1})\right)
\eta\left(J_2(\nu_1)\right)
\]
does not belong to $GE_2(K[z_1,z_2])$.

\begin{theorem}\label{Anick automorphism is wild}
{\rm (Umirbaev \cite{U2})}
The Anick automorphism
\[
\nu_1=(x+z(xz-zy),y+(xz-zy)z,z)
\]
is a wild automorphism of the free metabelian algebra $C$
and hence of the free associative algebra $K\langle x,y,z\rangle$.
\end{theorem}

In the same way Umirbaev showed the following theorem which
replaces the $z$-wildness of Theorem \ref{linear z-wild automorphisms}
with wildness in the usual sense.

\begin{theorem}\label{wild linear z-automorphisms by Umirbaev}
A linear $z$-automorphism of $K\langle x,y,z\rangle$
is tame as an automorphism of $K\langle x,y,z\rangle$
if and only if its $z$-Jacobian matrix belongs to $GE_2(K[z_1,z_2])$.
\end{theorem}

As in the case of polynomial algebras, the next step of studying automorphisms of free algebras
is to study coordinates, or primitive elements.
In particular, the following problem is motivated by \cite{UY}.
We state it in the case of three variables only, see \cite{DY4, DY5} for the general case.

\begin{problem}\label{wild coordinates}
If $f(x,y,z)\in K\langle x,y,z\rangle$ is an image of $x$ under a
$z$-wild automorphism, is there a tame automorphism (maybe not fixing $z$)
which also sends $x$ to $f(x,y,z)$?
\end{problem}

If such a tame automorphism does not exist, we call $f(x,y,z)$
a wild coordinate of $K\langle x,y,z\rangle$.

As an analog of the Strong Nagata Conjecture in \cite{UY}, we state

\begin{conjecture} {\bf (Strong Anick Conjecture)}
There exist wild coordinates in $K\langle x,y,z\rangle$.
In particular, the two nontrivial coordinates of the Anick automorphism
are both wild.
\end{conjecture}

The following theorem is one of the main results in Drensky and Yu \cite{DY4, DY5}.

\begin{theorem}\label{general form of strong Anick conjecture}
Let $K$ be a field of characteristic $0$ and let
the polynomial $f(x,y,z)\in K\langle x,y,z\rangle$ be linear in $x,y$.
If there exists a wild automorphism of $K\langle x,y,z\rangle$ which
fixes $z$
and sends $x$ to $f(x,y,z)$, then every automorphism of
$K\langle x,y,z\rangle$ which
sends $x$ to $f(x,y,z)$ is also wild.
So, $f(x,y,z)$ is a wild coordinate of $K\langle x,y,z\rangle$.
\end{theorem}

Theorem \ref{general form of strong Anick conjecture}
gives an algorithm  deciding whether a polynomial
$f(x,y,z)\in K\langle x,y,z\rangle$
which is linear in $x$ and $y$, is a tame coordinate. If it is, then
the algorithm finds a product of $z$-elementary automorphisms which sends $x$ to
$f(x,y,z)$. (Of course, as in the other algorithmic considerations we assume that the ground
field $K$ is
constructive, and we may perform calculations there.)
For this purpose, it is sufficient to combine
Theorem \ref{general form of strong Anick conjecture}
with the following convenient form of the result in
\cite{P, THK}, as stated in \cite{THK}.
Let $a(z_1,z_2),b(z_1,z_2)$ be two polynomials in $K[z_1,z_2]$. Then there exist
$c(z_1,z_2),d(z_1,z_2)\in K[z_1,z_2]$
such that the matrix
\[
G=\left(\begin{matrix}
a&c\\
b&d\\
\end{matrix}\right)
\]
belongs to $GE_2(K[z_1,z_2])$ if and only if we can bring the pair
$(a,b)$ to $(\alpha,0)$, $\alpha\in K^{\ast}$, using the Euclidean
algorithm only.

The following consequence of Theorem
\ref{general form of strong Anick conjecture} gives the affirmative
answer to the
Strong Anick Conjecture.

\begin{theorem}{\rm (Drensky and Yu \cite{DY4, DY5})}
The Strong Anick Conjecture is true. Namely, there exist wild
coordinates
in $K\langle x,y,z\rangle$. In particular,
the two nontrivial coordinates $x+z(xz-zy)$ and $y+(xz-zy)z$ of the
Anick automorphism
\[
\omega=(x+z(xz-zy),y+(xz-zy)z,z)
\]
are both wild.
\end{theorem}

We call an automorphism $\varphi=(f(x,y,z),g(x,y,z),z)$
of $K\langle x,y,z\rangle$
{\it Anick-like} if $f(x,y,z)$ and $g(x,y,z)$ are linear in $x,y$ and
the matrix
$J_z(\varphi)$ does not belong to $GE_2(K[z_1,z_2])$.
The following corollary from \cite{DY4, DY5}
is an analogue of a result from \cite{UY}.

\begin{corollary}
The two nontrivial coordinates $f(x,y,z),g(x,y,z)$ of any Anick-like
automorphism
\[
\varphi=(f(x,y,z),g(x,y,z),z)
\]
of $K\langle x,y,z\rangle$ are wild.
\end{corollary}

In the spirit  of the above results, we obtain
in \cite{DY4, DY5} the following theorem which is much stronger.

\begin{theorem}\label{more general theorem}
Let $f(x,y,z)$ be a $z$-coordinate of $K\langle x,y,z\rangle$
without terms depending only on $z$
(i.e., $f(0,0,z)=0$). If the linear part (with respect to $x$ and $y$)
$f_1(x,y,z)$ of $f(x,y,z)$
 is a $z$-wild coordinate,
then $f(x,y,z)$ itself is also a wild coordinate of $K\langle
x,y,z\rangle$.
\end{theorem}

The restriction $f(0,0,z)=0$ is essential for the proof of
Theorem \ref{more general theorem}.
Nevertheless, it seems very unlikely to
have a wild automorphism
$(f,g,z)$ with $f(0,0,z)=0$ such that $f(x,y,z)+a(z)$ is a
tame coordinate
for some polynomial $a(z)$ in view of the next theorem.

\begin{theorem}\label{nonlinear automorphisms with z fixed}
{\rm (Drensky and Yu \cite{DY4, DY5})}
Let $(f,g,z)$ be an automorphism of $K\langle x,y,z\rangle$
and let the linear
part (with respect to $x$ and $y$) of it, $(f_1,g_1,z)$,  be a $z$-wild
automorphism.
Then $(f,g,z)$ is also a wild automorphism of $K\langle x,y,z\rangle$.
\end{theorem}

The above theorem is much stronger
than the main result in \cite{U2} where only the
automorphisms linear with respect to $x$  and $y$ are dealt.

An important partial case of Example \ref{nonmetabelian example generalizing Anick}
gives a large class of wild automorphisms
and wild coordinates. (Recall that Corollary \ref{nonmetabelian generalization of Anick}
gives only the $z$-wildness of the automorphisms.)

\begin{example}\label{example generalizing Anick}
{\rm (Drensky and Yu \cite{DY4, DY5})}
Let $h(t,z)\in K\langle t,z\rangle$ and let $h(0,0)=0$.
If the linear component (with respect to $x,y$)
$h_1(xz-zy,z)$ of $h(xz-zy,z)$ is not equal to 0,
then the automorphism
\[
\sigma_h=(x+zh(xz-zy,z),y+h(xz-zy,z)z,z)
\]
belongs to the class of wild automorphisms
in Theorem \ref{more general theorem}.
\end{example}

\begin{example}\label{example without linear terms}
{\rm (Drensky and Yu \cite{DY4, DY5})}
A minor modification of the Anick automorphism is the automorphism
of $K\langle x,y,z\rangle$
\[
\omega_m=(x+z(xz-zy)^m,y+(xz-zy)^mz,z).
\]
Note that the automorphisms $\omega_m$, $m>1$,
are not covered by Theorem \ref{more general theorem}, as
the polynomials $z(xz-zy)^m$ and $(xz-zy)^mz$ have no
linear components with respect to $x$ and $y$.
\end{example}

\begin{theorem}\label{generalization of Anick automorphism}
{\rm (Drensky and Yu \cite{DY4, DY5})}
The above automorphisms $\omega_m$ are wild for all $m\geq 1$.
\end{theorem}

\begin{problem}
Are the two nontrivial coordinates of the above automorphism
$\omega_m$, $m>1$, both wild?
\end{problem}

The most general form of the result of
Umirbaev \cite{U2} gives that the automorphism $\vartheta=(f,g,h)$ of
the free metabelian algebra $M(x,y,z)$ is wild, if
it induces the identity automorphism of $K[x,y,z]$
and the matrix $J_2(\vartheta)(z_1,z_2)$
cannot be presented as a product of elementary matrices with entries
from $K[z_1,z_2]$, see Theorem \ref{criterion for tameness in metabelian case}.
Hence the classes of wild automorphisms and wild coordinates in
Theorem \ref{more general theorem},
Example \ref{example generalizing Anick}
and Example \ref{example without linear  terms} are
not covered by Umirbaev \cite{U2}.

Now we are going to show that at least two coordinates of the
automorphisms
of the class of Umirbaev are wild.

\begin{theorem}\label{the case of the theorem of Umirbaev}
{\rm (Drensky and Yu \cite{DY4, DY5})}
Let $\vartheta=(f,g,h)$ be an automorphism of
the free metabelian algebra $M(x,y,z)$ which
induces the identity automorphism of $K[x,y,z]$
and the matrix $J_2(\vartheta)(z_1,z_2)$ does not belong to
$GE_2(K[z_1,z_2])$.
Then the two  coordinates $f(x,y,z)$ and $g(x,y,z)$ are both wild.
\end{theorem}

The above results suggest the following problems.

\begin{problem}
Is it true that the two nontrivial coordinates of
a wild automorphism of $K\langle x,y,z\rangle$  fixing $z$ are both
wild?
\end{problem}

\begin{problem}
Is it true that every wild automorphism of $K\langle x,y,z\rangle$
contains at least two wild coordinates?
\end{problem}

Finally, we want to mention that the free metabelian algebra $C$
has quite peculiar wild automorphisms which cannot be lifted to automorphisms
of $K\langle x,y,z\rangle$. For example, such an automorphism is
\[
\rho=(x+x^2[y,z],y,z).
\]
Note that its Jacobian matrix
\[
J_M(\tau)=\left(\begin{matrix}
1&0&0\\
x_1^2(z_2-z_1)&1&0\\
x_1^2(y_1-y_2)&0&1\\
\end{matrix}\right)
\]
is a product of two elementary matrices. We refer to \cite{DY5} for details.

\begin{problem}
Is the polynomial $x+x^2[y,z]$ a wild coordinate of $C$?
Can it be lifted to a coordinate of $K\langle x,y,z\rangle$?
\end{problem}

\begin{problem}
Do there exist wild automorphisms and wild coordinates
of the free metabelian algebra $C_n$
of rank $n>3$? Are there wild automorphisms similar to the automorphism
$\rho$ constructed above?
\end{problem}

\section{Open Problems and Conjectures}

For some of the problems below we need $\text{char }K=0$.
For simplicity of the exposition, we consider only
the case of characteristic 0. Also, when we discuss algorithmic
problems, we assume that the base field is constructive. We fix
a finite set $X=\{x_1,\ldots,x_n\}$.
Some of the problems, considered for $X=\{x,y,z\}$
have two versions, depending whether we assume that
the automorphisms fix $z$ (or the derivations
annihilate $z$).

Yagzhev \cite{Y} suggested an algorithm which determines whether
an endomorphism of $K\langle X\rangle$ is an automorphism.
Hence we can check which $z$-endomorphisms of $K\langle x,y,z\rangle$
are automorphisms. Also, we have seen how to recognize the $z$-tame
automorphisms and coordinates
of $K\langle x,y,z\rangle$.
In the view of the results in \cite{SY1, DY2,SU1, SU3} exposed in the
survey part of the present paper, we suggest:

\begin{problem}
Find algorithms which decide whether a given polynomial in $K\langle x,y,z\rangle$
is a coordinate for a $z$-wild automorphism.
\end{problem}

\begin{conjecture}
If $\varphi$ is a wild $z$-automorphism of $K[x,y,z]$, then
it cannot be lifted to a $z$-automorphism (to any automorphism) of $K\langle x,y,z\rangle$.
\end{conjecture}

By the theorem of Umirbaev \cite{U1}, every automorphism of the polynomial
algebra $K[X]$ can be lifted to an automorphism of the free metabelian
algebra $M(X)=K\langle X\rangle/C^2$, where $C$ is the commutator ideal of
$K\langle X\rangle$. A concrete lifting (which also preserves $z$)
of the Nagata automorphism from $K[x,y,z]$ to $M(x,y,z)$ was given by the authors
and Gutierrez \cite{DGY}. It follows easily from the theory of finitely generated
PI-algebras that the theorem of Umirbaev holds also for
$K\langle X\rangle/T(R)$, where $T(R)$ is the T-ideal of the polynomial identities
(depending on the finite number of variables $X$) of a PI-algebra $R$
which satisfies a nonmatrix polynomial identity. The same lifting of the Nagata
automorphism to $M(x,y,z)$ can be used as a lifting to
$K\langle x,y,z\rangle/T(R)$.

\begin{problem} Describe the ($z$-)automorphisms
of $K[x,y,z]$ which can be lifted to automorphisms of
$K\langle x,y,z\rangle/T(M_n(K))$ and $K\langle x,y,z\rangle/T^2(M_n(K))$,
where $T(M_n(K))$ is the set of all polynomial identities of the
$n\times n$ matrix algebra. If all ($z$-)automorphisms
of $K[x,y,z]$ can be lifted to
$K\langle x,y,z\rangle/T^2(M_n(K))$, then they can be lifted also to
($z$-)automorphisms of any relatively free algebra
$K\langle x,y,z\rangle/T(R)$.
The case $K\langle x,y,z\rangle/T(M_2(K))$
is especially interesting because the structure of the algebra
$K\langle X\rangle/T(M_2(K))$ is well understood and one knows
how to work effectively there.
\end{problem}

In the case of the polynomial algebra $K[X]$,
if $\delta$ is a locally nilpotent derivation and
$K[X]^{\delta}$ is the algebra of constants (i.e., the
kernel) of $\delta$, then for any
$w\in K[X]^{\delta}$, the mapping $\Delta=w\delta$
is again a locally nilpotent derivation of $K[X]$
with the same algebra of constants as $\delta$. This fact is used
to construct new automorphisms of $K[X]$, including the Nagata
automorphism and to prove that some classes of automorphisms
are stably tame, see the paper by Smith \cite{Sm},
for better understanding of the automorphisms of $K[X]$, see
the book by van den Essen \cite{E2}, in the proof of Makar-Limanov of
the tameness of the automorphisms of $K[x,y]$, see \cite{ML3}
or the book by one of the authors \cite{D}, etc.
This construction does not produce anymore locally nilpotent derivations
in $K\langle X\rangle$ because for the proof that $\Delta=w\delta$
is a derivation one needs that $w$ is in the centre of the algebra.

\begin{problem} Starting from a locally nilpotent derivation
$\delta$ of $K\langle X\rangle$, find constructions which
produce other locally nilpotent derivations, preferable with the
same kernel as $\delta$. Can every locally nilpotent derivation
of $K[X]$ be lifted to a locally nilpotent derivation of
$K\langle X\rangle$?
\end{problem}

\begin{problem}
Can every locally nilpotent derivation
of $K[X]$ be lifted to a locally nilpotent derivation of
$K\langle X\rangle/T(R)$, where $R$ is a PI-algebra?
Can the derivation $\Delta=(y^2+xz)\delta$,
$\delta:(x,y,z)\to (-2y,z,0)$, (which is crucial
for the stable tameness of the Nagata automorphism) be lifted
to a locally nilpotent derivation of the free metabelian algebra $M_3(x,y,z)$?
\end{problem}

Recently van den Essen and Peretz \cite{EP} studied the automorphisms
$\varphi$ of $K[X]$ such that $\varphi(x_j)=\sum_{i=1}^nw_{ij}x_i$,
where each $w_{ij}$ is $\varphi$-invariant, i.e.,
these automorphisms are ``linear'' with coefficients from the algebra
$K[X]^{\varphi}$ of fixed points. They related them with locally
nilpotent derivations.

\begin{problem}
Describe, in the spirit of \cite{EP},
the automorphisms $\varphi$ of $K\langle X\rangle$
with the property that
\[
\varphi(x_j)=\sum_{i=1}^n\sum_pa_{ijp}x_ib_{ijp},\quad
a_{ijp},b_{ijp}\in K\langle X\rangle^{\varphi}.
\]
\end{problem}

In the commutative case Lamy \cite{L1, L2} described the ``nonaffine
orthogonal group'', i.e., the group of automorphisms
of $K[x,y,z]$ which fixes a nondegenerate quadratic form.
Considering the form $y^2+xz$, he proved that the group contains
automorphisms which do not belong to the subgroup
generated by Nagata-like automorphisms. Recall also the
general construction due to Freudenburg \cite{F1, F2} of locally nilpotent derivations
of $K[x,y,z]$ which do not annulate any coordinate.
Among them are derivations annulating the form $y^2+xz$.

\begin{problem} Describe the automorphisms of $K\langle x,y,z\rangle$
which fix the form $2y^2+xz+zx$. Describe the automorphisms with
this property which fix also $z$.
\end{problem}

\end{document}